\documentclass[11pt]{article}
\usepackage[T1]{fontenc}
\usepackage{lmodern}
\usepackage[colorlinks=true,linkcolor=RoyalBlue,urlcolor=RoyalBlue,citecolor=red]{hyperref}
\usepackage{amsmath,amsthm,amssymb,amsfonts}
\usepackage{enumerate}
\usepackage{epsfig}
\usepackage[dvipsnames,svgnames,table]{xcolor}
\usepackage{pgf,tikz}
\usetikzlibrary{calc}
\usepackage{fullpage}

\usepackage{booktabs}
\usepackage[nameinlink]{cleveref}
\usepackage{enumitem}
\usepackage{blkarray}
\usepackage{tikz-cd}

\usepackage{thmtools, thm-restate}

\newtheorem{theorem}{Theorem}[section]
\newtheorem{proposition}[theorem]{Proposition}
\newtheorem{lemma}[theorem]{Lemma}
\newtheorem{corollary}[theorem]{Corollary}

\newtheorem{claim}[theorem]{Claim}

\theoremstyle{definition}
\newtheorem{definition}[theorem]{Definition}

\theoremstyle{remark}
\newtheorem{remark}[theorem]{Remark}

\newenvironment{claimproof}[1][\proofname]
	{
		\proof[#1]%
			
	}
	{
		\endproof
	}

\newcommand{\CC}{\mathbb{C}}
\newcommand{\RR}{\mathbb{R}}

\newcommand{\FF}{\mathbb{F}}

\newcommand{\cR}{\mathcal{R}}

\DeclareMathOperator{\cl}{cl}

\tikzstyle{vertex}=[fill=black,circle,inner sep=0pt, minimum size=4pt]
\tikzstyle{edge}=[line width=1.5pt,black]

\title{Odd realisation numbers of rigid graphs}
\author{Sean Dewar\thanks{Department of Computer Science, KU Leuven. E-mail: \texttt{sean.dewar@kuleuven.be}}  \and Anthony Nixon\thanks{School of Mathematical Sciences, Lancaster University. E-mail: \texttt{a.nixon@lancaster.ac.uk}}  \and Ben Smith\thanks{School of Mathematical Sciences, Lancaster University. E-mail: \texttt{b.smith9@lancaster.ac.uk}}}

\begin{document}
\date{}
\maketitle

\begin{abstract}
Determining the number of realisations of a rigid graph for a specific choice of edge lengths is a fundamental problem in discrete geometry.
Gortler, Healy and Thurston showed that the number of real realisations of a generic framework is always even or equal to one.
Despite being rare, there do exist graphs whose number of complex realisations is odd and not equal to one.
In this article, we exhibit two necessary conditions for a rigid graph to have an odd realisation number other than one.
We also provide an infinite family of rigid graphs with an odd realisation number other than one.
\end{abstract}



\section{Introduction}\label{sec:intro}

A (bar-joint) framework $(G,p)$ in $\mathbb{R}^d$ is the combination of a finite, simple graph $G=(V,E)$ and a realisation $p:V\rightarrow \mathbb{R}^d$.
The framework is \emph{rigid} if the only edge-length-preserving continuous deformation of the vertices arises from isometries of $\mathbb{R}^d$. 
Determining the rigidity of frameworks is a computationally difficult problem \cite{Abbot}, and so one often assumes a genericity hypothesis on $p$. 
Asimow and Roth \cite{AsimowRothI} showed that rigidity of a graph is a generic property: $(G,p)$ is rigid for some generic realisation $p$ if and only if it is rigid for \emph{every} generic realisation.
As such, we say a graph $G$ is \emph{$d$-rigid} if there exists a generic rigid framework $(G,p)$ in $\mathbb{R}^d$.

A natural refinement of this question is the following: given a framework $(G,p)$, how many edge-length equivalent frameworks $(G,q)$ exist up to isometries of $\mathbb{R}^d$?
We denote this value by $r_d(G,p)$:
if $p$ is generic, then it follows that $(G,p)$ is rigid if and only if $r_d(G,p)$ is finite.
Many applications desire frameworks with $r_d(G,p)=1$, i.e. the framework has a unique realisation up to isometry; see for example \cite{SingerCucuringu:2010,Sensors}.
The combined results of Connelly \cite{conggr} and Gortler, Healy and Thurston \cite{gortler2010characterizing} give that unique realisability is also a generic property.
We thus say that $G$ is \emph{globally $d$-rigid} if there exists a generic realisation $p$ such that $r_d(G,p) = 1$.
As global rigidity cannot often be guaranteed in applications, $r_d(G,p)$ is a useful measure of how far a framework is from global rigidity.

\begin{figure}[tp]
	\begin{center}
        \begin{tikzpicture}[scale=0.8]
        \begin{scope}[xshift=0]
			\node[vertex] (1) at (0,1) {};
			\node[vertex] (2) at (0,-1) {};
			
			\node[vertex] (3) at (-1.5,-0.5) {};
			\node[vertex] (4) at (1.5,0.5) {};
			\node[vertex] (5) at (-0.5, 1.5) {};
			
			\draw[edge] (1)edge(2);
			\draw[edge] (1)edge(3);
			\draw[edge] (1)edge(4);
			\draw[edge] (2)edge(3);
			\draw[edge] (2)edge(4);
			\draw[edge] (5)edge(3);
			\draw[edge] (5)edge(4);
        \end{scope}
		\begin{scope}[xshift=100]
			\node[vertex] (1) at (0,1) {};
			\node[vertex] (2) at (0,-1) {};
			
			\node[vertex] (3) at (1.5,-0.5) {};
			\node[vertex] (4) at (1.5,0.5) {};
			\node[vertex] (5) at (3.68,0) {};
			
			\draw[edge] (1)edge(2);
			\draw[edge] (1)edge(3);
			\draw[edge] (1)edge(4);
			\draw[edge] (2)edge(3);
			\draw[edge] (2)edge(4);
			\draw[edge] (5)edge(3);
			\draw[edge] (5)edge(4);
		\end{scope}
        \begin{scope}[xshift=275]
            \draw[line width=1pt] (0,1.5) -- (0,-1.5);
		\end{scope}
        \begin{scope}[xshift=350]
			\node[vertex] (1) at (0,1) {};
			\node[vertex] (2) at (0,-1) {};
			
			\node[vertex] (3) at (1.5,-0.5) {};
			\node[vertex] (4) at (1.5,0.5) {};
			\node[vertex] (5) at (2,0) {};
			
			\draw[edge] (1)edge(2);
			\draw[edge] (1)edge(3);
			\draw[edge] (1)edge(4);
			\draw[edge] (2)edge(3);
			\draw[edge] (2)edge(4);
			\draw[edge] (5)edge(3);
			\draw[edge] (5)edge(4);
        \end{scope}
		\end{tikzpicture}
	\end{center}
    \caption{Two different choices of edge lengths for the same graph. The left choice gives $r_2(G,p) = 4$: the two realisations pictured plus two more via reflecting the degree 2 vertex through the line passing through its neighbours. The right choice only gives $r_2(G,p)=2$: again, the additional realisation can be found via reflecting the degree 2 vertex.}\label{fig:notgeneric}
\end{figure}
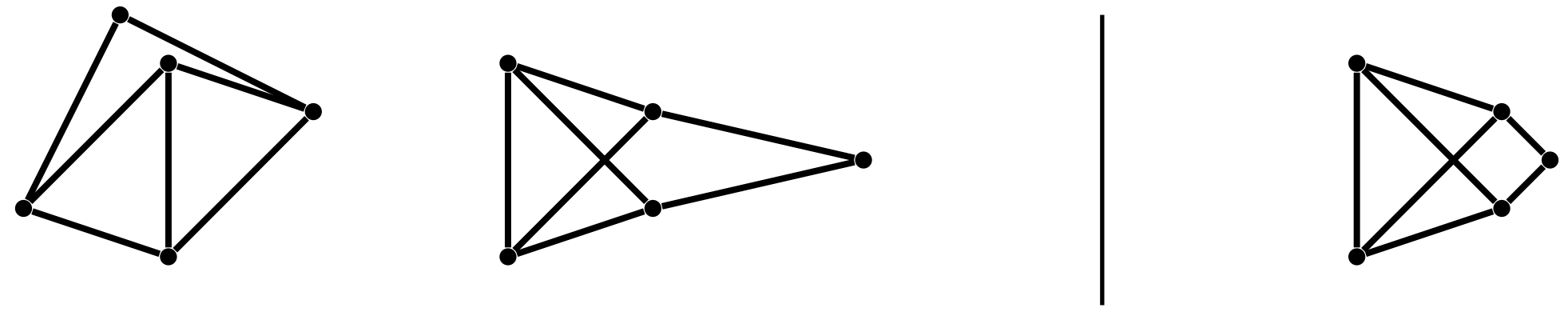

Unlike rigidity and global rigidity, the value $r_d(G,p)$ is not a generic property: see \Cref{fig:notgeneric} for an example.
There are two alternative ways to associate a single `generic realisation number' to a graph $G$:
\begin{enumerate}
    \item We define the \emph{real $d$-realisation number} $r_d(G)$ to be the largest value of $r_d(G,p)$ over all generic realisations $p \colon V \rightarrow \RR^d$.
    \item We extend $r_d(G,p)$ to additionally count complex realisations $p \colon V \rightarrow \CC^d$. This value \emph{is} a generic property, which we define to be the \emph{(complex) $d$-realisation number} $c_d(G)$.
\end{enumerate}
 It is immediate that $c_d(G) \geq r_d(G)$ for all graphs.
In general, both real and complex realisation numbers are often difficult to compute, and there have been a myriad of techniques developed in the last two decades to provide exact values and bounds for both $c_d(G)$ and $r_d(G)$:
see for example \cite{BARTZOS2021189,BartzosUpper,cggkls,JO19,Stef10}.

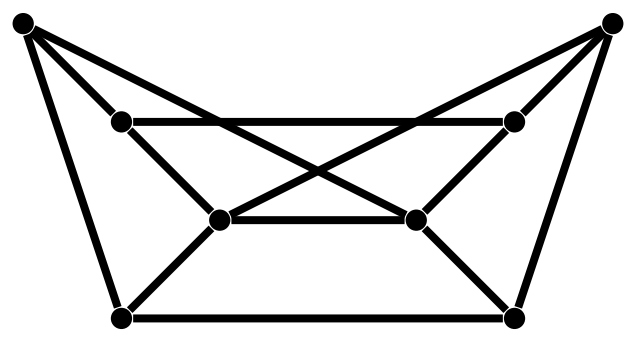
\begin{figure}[tp]
	\begin{center}
            \begin{tikzpicture}[scale=0.65]
			\node[vertex] (11) at (-1,0) {};
			\node[vertex] (21) at (-2,1) {};
			\node[vertex] (31) at (-2,-1) {};
			
			\node[vertex] (12) at (1,0) {};
			\node[vertex] (22) at (2,1) {};
			\node[vertex] (32) at (2,-1) {};

			\node[vertex] (a) at (-3,2) {};
                \node[vertex] (b) at (3,2) {};
			
			\draw[edge] (11)edge(12);
			\draw[edge] (21)edge(22);
			\draw[edge] (31)edge(32);
			
			\draw[edge] (11)edge(21);
			\draw[edge] (31)edge(11);
			
			\draw[edge] (12)edge(22);
			\draw[edge] (32)edge(12);

			\draw[edge] (a)edge(21);
			\draw[edge] (a)edge(12);
			\draw[edge] (a)edge(31);

                \draw[edge] (b)edge(22);
			\draw[edge] (b)edge(11);
			\draw[edge] (b)edge(32);
		\end{tikzpicture}
	\end{center}
	\caption{A graph $G$ with $c_2(G)=45$. This is the smallest $2$-rigid graph with $c_2(G)$ odd and greater than 1.}\label{fig:smallestodd}
\end{figure}

An interesting difference between the real and complex realisation numbers arises in their possible parity.
In \cite{conggr}, Connelly introduced the full rank stress sufficient condition for generic framework global rigidity, which was later proven to be necessary by Gortler, Healy and Thurston \cite{gortler2010characterizing}.
The latter's proof relies on proving that, if a generic framework does not have a full rank stress, then it has an even number of equivalent realisations up to congruence.
This has an interesting consequence for the real realisation number: for any $d$-rigid graph $G$, either $r_d(G)$ is even or is equal to one. 
From this, it would be natural to assume that the complex $d$-realisation number of any graph that is not globally $d$-rigid would also be even.
This is, however, not true: the graph pictured in \Cref{fig:smallestodd} is minimally 2-rigid with $c_2(G) = 45$.

Despite this example, there is much evidence to suggest that vanishingly few graphs have odd complex realisation number other than one.
When $d=1$, the realisation number has a straightforward combinatorial formula that implies $c_1(G)$ is always even when $G$ is not globally $1$-rigid.
When $d=2$, the 2-realisation number $c_2(G)$ has been calculated for all minimally 2-rigid graphs on at most twelve vertices \cite{capcodata}, and the table given in \Cref{fig:realisation+table} lists the number of minimally 2-rigid graphs and the number with $c_2$ odd. 
Heuristically, it appears that vanishingly-many have odd realisation number. 
In all dimensions $d$, probabilistic approaches have shown that asymptotically all $d$-rigid graphs have $c_d(G)$ equal to a power of two \cite{dewarnixonsmith_random}. 
Moreover, certain common extension moves that preserve rigidity are known to double the realisation number (see for example \Cref{lem:0ext}).
This suggests graphs with odd complex realisation number are exceedingly rare: a priori, there may only be finitely many.

\begin{table}[htp]
\centering
\begin{tabular}{ | c | c | c | c |}
  \hline
  $|V|$ & $\#$ min.~2-rigid & $\#$ $c_2$ odd  & $\%$ odd \\ 
  \hline
  7 & 70 & 0 & 0$\%$ \\ 
  \hline
  8 & 608 & 1 & 0.164$\%$ \\ 
  \hline
  9 & 7222 & 3 & 0.042$\%$ \\ 
  \hline
  10 & 110132 & 13 & 0.012$\%$ \\ 
  \hline
  11 & 2039273 & 153 & 0.0075$\%$ \\ 
  \hline
  12 & 44176717 & 2077 & 0.0047$\%$ \\ 
  \hline
\end{tabular}
\caption{A table of minimally 2-rigid graphs and the proportion with odd realisation number.}
\label{fig:realisation+table}
\end{table}

\subsection{Our contributions}

We exhibit two different structural properties that are necessary for graphs with odd realisation number.
Inspired by \cite{JO19}, we say a graph $G$ is \emph{$d$-Minkowski} if for every real generic framework $(G,p)$, there exists a non-real framework $(G,q)$ with the same edge lengths whose coordinates $q_i(v)$ are either real or purely imaginary.
We show that any $d$-rigid graph whose $d$-realisation number is odd and not equal to one is necessarily $d$-Minkowski (\Cref{prop:mink1}).

For our second necessary condition, we generalise the following definition of Penne \cite{Penne}.
We say that a $d$-rigid graph is a \emph{minimally $d$-isostatic graph} ($d$-MIG) if its only proper $d$-rigid subgraphs are complete graphs on at most $d$ vertices.
For example, $K_{3,3}$ is a 2-MIG as it is $2$-rigid and its only $2$-rigid subgraphs are single edges.
We prove that being a $d$-MIG is a necessary condition for minimally $d$-rigid graphs with odd $d$-realisation number other than one.

\begin{restatable}{theorem}{mig}\label{thm:mig}
    Let $G$ be a minimally $d$-rigid graph with at least $d+2$ vertices where $c_d(G)$ is odd.
    Then $G$ is a $d$-MIG.
\end{restatable}

We improve upon this result in the special case of $d=2$ by removing the minimality assumption from this statement (\Cref{cor:mig2d}). That is, any 2-rigid graph with odd realisation number $>1$ is actually minimally 2-rigid and in fact a 2-MIG.

It follows from \cite{dewarnixonsmith_random} that $d$-rigid but not globally $d$-rigid graphs with an odd $d$-realisation number are rare.
In fact, it is unclear if there exist infinitely many for any particular choice of $d$: indeed, if $d=1$, then there exist no such graphs.
We show that there are infinitely many such graphs when $d=2$.

\begin{restatable}{theorem}{infinite}\label{thm:allodd}
    There are infinitely many minimally 2-rigid graphs with an odd 2-realisation number.
\end{restatable}

\begin{figure}[tp]
	\begin{center}
            \begin{tikzpicture}
            \begin{scope}[scale=0.75, xshift = -100]
			\node[vertex] (00) at (0,0) {};
			\node[vertex] (01) at (0,1) {};
			\node[vertex] (02) at (0,2) {};
                
			\node[vertex] (10) at (1,0) {};
			\node[vertex] (11) at (1,1) {};
			\node[vertex] (12) at (1,2) {};

			\draw[edge] (00)edge(10);

			\draw[edge] (01)edge(11);

			\draw[edge] (02)edge(12);

			\draw[edge] (00)edge(01);
			\draw[edge] (10)edge(11);

			\draw[edge] (01)edge(12);
			\draw[edge] (11)edge(02);

            \draw[edge] (00) to [bend left = 45] (02);

            \draw[edge] (10) to [bend right = 45] (12);
            \end{scope}

            \begin{scope}[scale=0.75]
			\node[vertex] (00) at (0,0) {};
			\node[vertex] (01) at (0,1) {};
			\node[vertex] (02) at (0,2) {};
                
			\node[vertex] (10) at (1,0) {};
			\node[vertex] (11) at (1,1) {};
			\node[vertex] (12) at (1,2) {};

			\node[vertex] (20) at (2,0) {};
			\node[vertex] (21) at (2,1) {};
			\node[vertex] (22) at (2,2) {};

			\draw[edge] (00)edge(10);
			\draw[edge] (10)edge(20);

			\draw[edge] (01)edge(11);
			\draw[edge] (11)edge(21);

			\draw[edge] (02)edge(12);
			\draw[edge] (12)edge(22);

			\draw[edge] (00)edge(01);
			\draw[edge] (10)edge(11);
			\draw[edge] (20)edge(21);

			\draw[edge] (01)edge(12);
			\draw[edge] (11)edge(02);
            
			\draw[edge] (11)edge(22);
			\draw[edge] (21)edge(12);

            \draw[edge] (00) to [bend left = 45] (02);

            \draw[edge] (20) to [bend right = 45] (22);
            \end{scope}
            
            \begin{scope}[scale=0.75, xshift = 130]
			\node[vertex] (00) at (0,0) {};
			\node[vertex] (01) at (0,1) {};
			\node[vertex] (02) at (0,2) {};
			\node[vertex] (02) at (0,2) {};

			\node[vertex] (10) at (1,0) {};
			\node[vertex] (11) at (1,1) {};
			\node[vertex] (12) at (1,2) {};

			\node[vertex] (20) at (2,0) {};
			\node[vertex] (21) at (2,1) {};
			\node[vertex] (22) at (2,2) {};

                \node[vertex] (30) at (3,0) {};
			\node[vertex] (31) at (3,1) {};
			\node[vertex] (32) at (3,2) {};

			\draw[edge] (00)edge(10);
			\draw[edge] (10)edge(20);
			\draw[edge] (20)edge(30);

			\draw[edge] (01)edge(11);
			\draw[edge] (11)edge(21);
			\draw[edge] (21)edge(31);

			\draw[edge] (02)edge(12);
			\draw[edge] (12)edge(22);
			\draw[edge] (22)edge(32);

			\draw[edge] (00)edge(01);
			\draw[edge] (10)edge(11);
			\draw[edge] (20)edge(21);
			\draw[edge] (30)edge(31);

			\draw[edge] (01)edge(12);
			\draw[edge] (11)edge(02);
            
			\draw[edge] (11)edge(22);
			\draw[edge] (21)edge(12);

                \draw[edge] (21)edge(32);
			\draw[edge] (31)edge(22);

            \draw[edge] (00) to [bend left = 45] (02);

            \draw[edge] (30) to [bend right = 45] (32);
            \end{scope}

            \begin{scope}[scale=0.75, xshift = 285]
			\node[vertex] (00) at (0,0) {};
			\node[vertex] (01) at (0,1) {};
			\node[vertex] (02) at (0,2) {};
			\node[vertex] (02) at (0,2) {};

			\node[vertex] (10) at (1,0) {};
			\node[vertex] (11) at (1,1) {};
			\node[vertex] (12) at (1,2) {};

			\node[vertex] (20) at (2,0) {};
			\node[vertex] (21) at (2,1) {};
			\node[vertex] (22) at (2,2) {};

                \node[vertex] (30) at (3,0) {};
			\node[vertex] (31) at (3,1) {};
			\node[vertex] (32) at (3,2) {};

                \node[vertex] (40) at (4,0) {};
			\node[vertex] (41) at (4,1) {};
			\node[vertex] (42) at (4,2) {};

			\draw[edge] (00)edge(10);
			\draw[edge] (10)edge(20);
			\draw[edge] (20)edge(30);
			\draw[edge] (30)edge(40);

			\draw[edge] (01)edge(11);
			\draw[edge] (11)edge(21);
			\draw[edge] (21)edge(31);
			\draw[edge] (31)edge(41);

			\draw[edge] (02)edge(12);
			\draw[edge] (12)edge(22);
			\draw[edge] (22)edge(32);
			\draw[edge] (32)edge(42);

			\draw[edge] (00)edge(01);
			\draw[edge] (10)edge(11);
			\draw[edge] (20)edge(21);
			\draw[edge] (30)edge(31);
			\draw[edge] (40)edge(41);

			\draw[edge] (01)edge(12);
			\draw[edge] (11)edge(02);
            
			\draw[edge] (11)edge(22);
			\draw[edge] (21)edge(12);

                \draw[edge] (21)edge(32);
			\draw[edge] (31)edge(22);

                \draw[edge] (31)edge(42);
			\draw[edge] (41)edge(32);

            \draw[edge] (00) to [bend left = 45] (02);

            \draw[edge] (40) to [bend right = 45] (42);
            \end{scope}

                  \begin{scope}[scale=0.75, xshift = 100,
                  yshift = -100]
			\node[vertex] (00) at (0,0) {};
			\node[vertex] (01) at (0,1) {};
			\node[vertex] (02) at (0,2) {};

			\node[vertex] (10) at (1,0) {};
			\node[vertex] (11) at (1,1) {};
			\node[vertex] (12) at (1,2) {};

			\node[vertex] (20) at (2,0) {};
			\node[vertex] (21) at (2,1) {};
			\node[vertex] (22) at (2,2) {};

                \node[vertex] (30) at (3,0) {};
			\node[vertex] (31) at (3,1) {};
			\node[vertex] (32) at (3,2) {};

			\draw[edge] (00)edge(10);
			\draw[black,dashed] (10)--(20);
			\draw[edge] (20)edge(30);

			\draw[edge] (01)edge(11);
			\draw[black,dashed] (11)--(21);
			\draw[edge] (21)edge(31);

			\draw[edge] (02)edge(12);
			\draw[black,dashed] (12)--(22);
			\draw[edge] (22)edge(32);

			\draw[edge] (00)edge(01);
			\draw[edge] (10)edge(11);
			\draw[edge] (20)edge(21);
			\draw[edge] (30)edge(31);

			\draw[edge] (01)edge(12);
			\draw[edge] (11)edge(02);
            
			\draw[black,dashed] (11)--(22);
			\draw[black,dashed] (21)--(12);

                \draw[edge] (21)edge(32);
			\draw[edge] (31)edge(22);

            \draw[edge] (00) to [bend left = 45] (02);

            \draw[edge] (30) to [bend right = 45] (32);
            \end{scope} 
		\end{tikzpicture}
	\end{center}
	\caption{The family of minimally 2-rigid graphs $A_n$ for $n\geq 2$. 
    Above are $A_2, \dots, A_5$, with a general schematic of $A_n$ below.}
    \label{fig:odd_intro}
\end{figure}
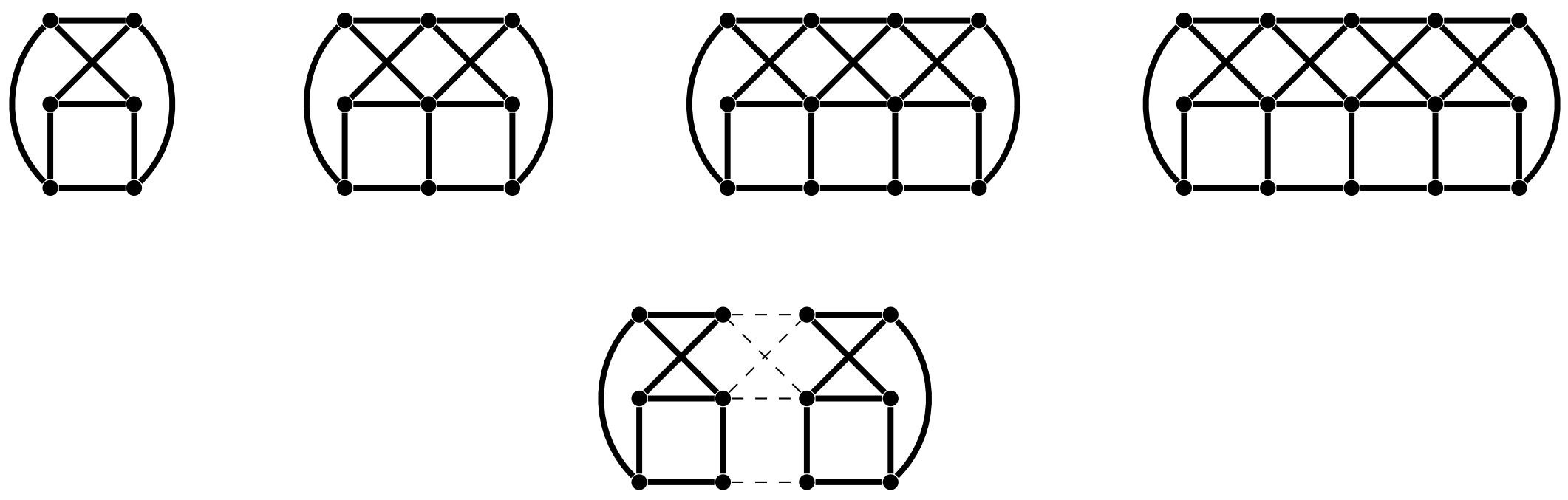

We give two complementary proofs of this theorem.
The first is constructive: we exhibit an infinite family of graphs $A_n$, depicted in \Cref{fig:odd_intro}, that are minimally 2-rigid and have $c_2(A_n)$ odd for all $n\geq 3$.
The second is non-constructive but more general: we introduce a gluing operation that takes as input two minimally 2-rigid graphs with odd 2-realisation number, and outputs two new minimally 2-rigid graphs, exactly one of which has odd 2-realisation number (\Cref{thm:oddglue}).
Iterating this gluing procedure provides a non-constructive proof of \Cref{thm:allodd} demonstrating a wider class of minimally 2-rigid graphs with odd 2-realisation number.
Both approaches rely on calligraph gluing operations introduced by Grasegger, El Hilany and Lubbes in \cite{ghl24}.

\subsection{Layout of paper}

The paper has the following layout.
We first review notation, terminology and provide preliminary results from rigidity theory and realisation numbers in \Cref{sec:prelims}.
In \Cref{sec:odd}, we give structural results on graphs with odd realisation numbers not equal to one.
In \Cref{sec:infiniteoddgraphs}, we utilise the concept of calligraphs to construct infinitely many minimally 2-rigid graphs with an odd 2-realisation number greater than 1.
In \Cref{sec:nonconstruct}, we introduce a gluing operation that produces two graphs, exactly one of which has an odd 2-realisation number greater than 1. 
This gluing operation provides a non-constructive proof that there exist infinitely many minimally 2-rigid graphs with an odd 2-realisation number greater than 1. 

\section{Preliminaries}
\label{sec:prelims}

We now fix $\mathbb{F}$ to be either the field of real numbers or the field of complex numbers.
If the specific choice of field is not stated then it is assumed to be $\mathbb{C}$.
We will equip each vector space $\mathbb{F}^d$ with the quadratic form $\| (x_1,\ldots,x_d) \|^2 := \sum_{i=1}^d x_i^2$.
The notation here is intentional, as when $\mathbb{F}=\mathbb{R}$ our chosen quadratic form is exactly the square of the standard Euclidean norm.
The isometries of the quadratic space $(\mathbb{F}^d,\|\cdot\|^2)$ are exactly the affine transformations $x \mapsto Ax +b$ where $b \in \mathbb{F}^d$ and $A \in O(d,\mathbb{F})$, i.e.~the group of $d \times d$ matrices with coefficients in $\mathbb{F}$ where $M^T M = M M^T = I_d$.

\subsection{Rigidity theory for real and complex frameworks}

Recall that a \emph{$d$-dimensional framework} $(G,p)$ in $\mathbb{F}^d$ is the combination of a finite, simple graph $G=(V,E)$ and a realisation $p:V\rightarrow \mathbb{F}^d$.
The linear space of all such realisations of $G$ is denoted by $(\mathbb{F}^d)^V$.
We say that two realisations $p,q$ of $G$ in $\mathbb{F}^d$ are \emph{equivalent} if $\|p(v) - p(w)\|^2 = \|q(v) - q(w)\|^2$ holds for every edge $vw \in E$.
Furthermore, we say that two realisations $p,q$ of $G$ in $\mathbb{F}^d$ are \emph{congruent} (denoted $p \sim q$) if there exists an isometry of $\mathbb{F}^d$ that sends $p$ to $q$.
For the case where $\mathbb{F}=\mathbb{R}$, congruence is equivalent to the condition that $\|p(v) - p(w)\|^2 = \|q(v) - q(w)\|^2$ holds for all $v, w \in V$.

A \emph{continuous flex} of a framework $(G,p)$ in $\mathbb{F}^d$ is a continuous path $\rho :(-\varepsilon,\varepsilon) \rightarrow (\mathbb{F}^d)^V$ where $\rho_0=p$ and each framework $(G,\rho_t)$ is equivalent to $(G,p)$. We say a continuous flex $\rho$ is \emph{trivial} if $\rho_t \sim p$ for all $t \in (-\varepsilon,\varepsilon)$.
We now say $(G,p)$ is \emph{rigid} if every continuous flex of $(G,p)$ is trivial.

We will be interested in the rigidity properties that hold for `almost all realisations' of $G$.
A realisation $p$ is said to be \emph{generic} if the set of coordinates of $p$ forms an algebraically independent set over $\mathbb{Q}$ of $d|V|$ elements.
Any realisation that is congruent to a generic realisation is said to be \emph{quasi-generic}.

Determining whether a real framework is rigid is coNP-hard, and so considered computationally difficult \cite{Abbot}.
However, when a framework $(G,p)$ is generic, we can reduce this problem to linear algebra by performing a standard linearisation technique on the length constraints.
More precisely, we define the \emph{rigidity matrix} $R(G,p)$ of a $d$-dimensional framework $(G,p)$ to be the $|E|\times d|V|$ matrix whose rows are indexed by the edges and $d$-tuples of columns indexed by the vertices. The row for an edge $e=uv$ is given by:
$$
\begin{bmatrix} 
   ~ 0 & \cdots & 0 & p(u)-p(v) & 0 & \cdots & 0 & p(v)-p(u) & 0 & \cdots & 0 ~
\end{bmatrix}
$$ 
where $p(u)-p(v)$ occurs over the $d$-tuple of columns indexed by $u$ and $p(v)-p(u)$ occurs in the $d$-tuple of columns indexed by $v$.
Maxwell \cite{Maxwell} observed that for real frameworks,
$\mbox{rank } R(G,p)\leq d|V|-\binom{d+1}{2}$ whenever $p$ affinely spans $\mathbb{R}^d$;
moreover, this observation can be extended to complex frameworks.
We say that $(G,p)$ is \emph{infinitesimally rigid} if $\mbox{rank } R(G,p) = d|V|-\binom{d+1}{2}$, or $(G,p)$ is a simplex (i.e.~$G=K_n$ for $n \leq d+1$) with affinely independent vertices.

A fundamental result of Asimow and Roth \cite{AsimowRothI} tells us that, when $p$ is generic and $G$ has at least $d+1$ vertices, $(G,p)$ is rigid in $\mathbb{R}^d$ if and only if it is infinitesimally rigid;
additionally, this result can be easily extended to complex frameworks.
The set of realisations in $(\mathbb{F}^d)^V$ which are not infinitesimally rigid frameworks is an algebraic set defined by rational coefficient polynomials, and hence both rigidity and infinitesimal rigidity are generic properties for frameworks in $\mathbb{F}^d$.
We observe here that a framework in $\mathbb{R}^d$ is infinitesimally rigid when considered as a framework in $\mathbb{R}^d$ if and only if it is infinitesimally rigid when considered as a framework in $\mathbb{C}^d$.
Hence, given generic frameworks $(G,p)$, $(G,q)$ in $\mathbb{R}^d$ and $\mathbb{C}^d$ respectively, we have that $(G,p)$ is (infinitesimally) rigid if and only if $(G,q)$ is (infinitesimally) rigid.
This motivates the following definition:

\begin{definition}
    We say that a graph $G=(V,E)$ is \emph{$d$-rigid} if there exists a generic framework $(G,p)$ in either $\mathbb{R}^d$ or $\mathbb{C}^d$ that is rigid.
    If, in addition, $G-e$ is not $d$-rigid for any edge $e \in E$, then we say that $G$ is \emph{minimally $d$-rigid}.
\end{definition}

A common technique for proving that graphs are $d$-rigid is constructing them via extension moves that preserve rigidity.
A graph $G'$ is said to be obtained from another graph $G$ by:
\begin{itemize}
    \item a \emph{($d$-dimensional) 0-extension} if $G=G'-v$ for a vertex  $v\in V(G')$ with $d_{G'}(v)=d$,
    \item a \emph{($d$-dimensional) 1-extension} if $G=G'-v+xy$ for a vertex $v\in V(G')$ with $d_{G'}(v)=d+1$ and $x,y\in N_{G'}(v)$ with $xy \notin E(G')$.
\end{itemize}
It is well known that both of these operations preserve minimal $d$-rigidity.

\begin{lemma}[see e.g., {\cite[Lem.~11.1.1, Thm 11.1.7]{Whiteley96}}]\label{lem:01ext}
Let $G$ be minimally $d$-rigid and let $G'$ be obtained from $G$ by a $d$-dimensional 0-extension or a 1-extension. 
Then $G'$ is minimally $d$-rigid.
\end{lemma}

It will be useful at points to utilise techniques from matroid theory to study generic rigidity; see \cite{Oxley} for an introduction.
The (generic) \emph{$d$-dimensional rigidity matroid} of a graph $G=(V,E)$ is the matroid $\mathcal{R}_d(G)$ on $E$ in which a set of edges $F\subseteq E$ is independent whenever the corresponding rows of $R(G,p)$ are independent, for some (or equivalently every) generic $p:V\to \mathbb{F}^d$. 
We will often just write $\mathcal{R}_d$ for the $d$-dimensional rigidity matroid on some suitably large complete graph. 
We denote the closure operator of $\mathcal{R}_d$ by $\cl_d$: recall that a graph $G = (V,E)$ is $d$-rigid if and only if $\cl_d(E) = E(K_V)$.
We say that a graph $G$ is \emph{$\mathcal{R}_d$-connected} if $\mathcal{R}_d(G)$ is a connected matroid, i.e. every pair of edges is contained in a common circuit of $\mathcal{R}_d(G)$.
We also say that a graph $G$ is a \emph{$d$-circuit} if its edge set forms a circuit of $\mathcal{R}_d(G)$.

In the case when $d=2$, the rigidity matroid $\cR_2$ has a simple combinatorial characterisation, which in turn gives combinatorial criteria for rigidity in 2-dimensions.

\begin{theorem}[\cite{PollaczekGeiringer,Laman,JacksonJordan:2005}]\label{thm:2rigmatroid}
    Let $G = (V,E)$ be a graph with $|V| \geq 2$.
    \begin{enumerate}
        \item \label{thm:2rigmatroid:ind} $E$ is independent in $\mathcal{R}_2$ if and only if $|E'| \leq 2|V'|-3$ for all subgraphs $G' = (V',E')$ of $G$ with at least two vertices.
        \item \label{thm:2rigmatroid:rigid} $G$ is minimally 2-rigid if and only if $E$ is independent in $\mathcal{R}_2$ and $|E| = 2|V|-3$.
        \item \label{thm:2rigmatroid:circ} If $G$ is $\mathcal{R}_2$-connected, then $G-e$ is 2-rigid for any choice of edge $e \in E$.
    \end{enumerate}
\end{theorem}

The last bit of terminology we require from classical rigidity theory is the following.
A framework $(G,p)$ in $\FF^d$ is \emph{globally rigid} if every equivalent framework to $(G,p)$ is congruent to $(G,p)$.
It was proven by Connelly \cite{conggr} and Gortler, Healy and Thurston \cite{gortler2010characterizing} that a single generic framework $(G,p)$ in $\RR^d$ is globally rigid if and only if every generic framework $(G,q)$ in $\RR^d$ is globally rigid.
Gortler and Thurston \cite{Gortler2014} later extended this to frameworks in $\CC^d$, and further proved that a generic globally rigid framework in $\RR^d$ remains globally rigid when considered as a framework in $\CC^d$.
Since every generic framework in $\RR^d$ is also a generic framework in $\CC^d$, the following definition is consistent:
a graph $G$ is \emph{globally $d$-rigid} if and only if every (equivalently, some) generic framework $(G,p)$ in $\RR^d$ (equivalently, $\CC^d$) is globally rigid.

As with 2-rigidity, global 2-rigidity has a combinatorial characterisation.

\begin{theorem}[Jackson and Jord\'{a}n \cite{JacksonJordan:2005}]\label{thm:2globrig}
    A graph is globally 2-rigid if and only if it is a complete graph or it is both $\mathcal{R}_2$-connected and 3-connected.
\end{theorem}

An immediate corollary of \Cref{thm:2rigmatroid,thm:2globrig} is that minimally 2-rigid graphs on at least four vertices cannot be globally 2-rigid.

\subsection{Realisation numbers}

Given a graph $G = (V,E)$ that is $d$-rigid, we wish to study the number of equivalent realisations for a generic real framework $(G,p)$ up to isometry.
Formally, we define
\[
r_d(G,p) = \left|\left\{ q \in (\mathbb{R}^d)^V \; : \; \|p (v)-p(w)\|^2 = \|q (v)-q(w)\|^2 \text{ for all } vw \in E \right\}/ \sim \right| \, ,
\]
where $p \sim q$ if $p, q$ congruent.
As showcased in \Cref{fig:notgeneric}, this is no longer a generic property: $r_d(G,p)$ varies as we range over generic realisations $p \colon V\rightarrow \RR^d$.
As such, authors define the \emph{real $d$-realisation number} to be
\[
r_d(G) = \max \left\{r_d(G,p) \; : \; p\colon V \rightarrow \RR^d \text{ generic } \right\}\, .
\]
We make a point here that the combined results of Connelly \cite{conggr} and Gortler, Healy and Thurston \cite{gortler2010characterizing} imply that $r_d(G) = 1$ if and only if $r_d(G,p)=1$ for some generic $d$-dimensional framework $(G,p)$.

The value $r_d(G)$ is somewhat difficult to work with, as often there will exist generic realisations which do not have this number of equivalent frameworks.
To avoid this issue, we often instead approach this problem algebraically by counting complex realisations also.
We define
\begin{equation*}
	C_d(G,p) := \left\{ q \in (\mathbb{C}^d)^V \; : \; \|p (v)-p(w)\|^2 = \|q (v)-q(w)\|^2 \text{ for all } vw \in E\right\}/ \sim
\end{equation*}
to be the \emph{realisation space of $(G,p)$}.
The cardinality of $C_d(G,p)$ is constant over all generic frameworks \cite[Thm 2.4]{dewarsmithnixon2026}.

\begin{definition}\label{def:count}
	The \emph{(complex) $d$-realisation number} of a $d$-rigid graph $G=(V,E)$ is an element of $\mathbb{N} \cup \{\infty\}$ given by
	\begin{align*}
		c_d(G) :=		
		\begin{cases}
			|C_d(G,p)| \text{ for generic $p \in (\mathbb{C}^d)^V$} &\text{if } |V| \geq d+1, \\
			 1 &\text{if } |V| \leq d \text{ and $G$ is complete}, \\
			 \infty &\text{if } |V| \leq d \text{ and $G$ is not complete}.
		\end{cases}
	\end{align*}
\end{definition}

While it is immediate from the definitions that $c_d(G) \geq r_d(G)$, we note that these two quantities agree when $G$ is globally $d$-rigid.
Explicitly, Gortler and Thurston \cite{Gortler2014} showed that $G$ is globally $d$-rigid if and only if $c_d(G)=1$.

It will sometimes be beneficial to work with frameworks that have been `pinned' in a suitable position that kills almost all isometries beforehand.
Fixing a sequence of $d$ vertices $v_1,\ldots, v_d$ of $G$, we define
\begin{align}\label{eq:xset}
	X_{G,d} := \left\{ p \in (\mathbb{C}^{d})^V : \text{ for all } 1\leq k\leq d,~ p_j(v_k) =0 \text{ if } j \geq k  \right\}\,.
\end{align}
Intuitively, this is the space of frameworks where $v_1$ is pinned to the origin, $v_2$ is pinned to the $x$-axis, $v_3$ is pinned to the $xy$-plane, etc.

Almost every framework is congruent to a framework in $X_{G,d}$; in particular, every quasi-generic framework is congruent to a framework in $X_{G,d}$.
We note that the frameworks in $X_{G,d}$ are never generic as they have been pinned in special position, but $X_{G,d}$ does contain a dense subset of quasi-generic frameworks.

\begin{lemma}[\cite{dewarsmithnixon2026}]\label{lem:quasixgd}
    Let $(G,p)$ be a quasi-generic rigid framework in $\mathbb{C}^d$ with at least $d+1$ vertices.
    Then for any choice of vertices $v_1,\ldots,v_d$ defining the space $X_{G,d}$,
    there exists a congruent framework $(G,q)$ with $q \in X_{G,d}$.
\end{lemma}

Finally, we note that the $d$-realisation number doubles under 0-extensions.
This will allow us to easily compute the $d$-realisation number of certain classes of graphs.

\begin{lemma}[{\cite[Lem.~7.1]{dewar23}}]\label{lem:0ext}
    Let $G = (V,E)$ be a $d$-rigid graph with at least $d+1$ vertices.
    If $G'$ is obtained from $G$ via a ($d$-dimensional) 0-extension, then $c_d(G') = 2c_d(G)$.
\end{lemma}

\section{Structural properties of graphs with odd realisation numbers}\label{sec:odd}

\subsection{The Minkowski property}

For the following, $\mathfrak{i}$ denotes the imaginary number.

\begin{definition}\label{def:mink}
    A graph $G$ with at least $d+1$ vertices is \emph{$d$-Minkowski} if for every real generic $(G,p)$ in $\mathbb{R}^d$,
    there exists an equivalent 
    framework $(G,q)$ in $\mathbb{C}^d$ and a set $\emptyset \subsetneq I \subseteq [d]$ such that $q_i(v) \in \mathbb{R}$ if $i \in [d] \setminus I$ and $q_i(v) \in \mathfrak{i}\mathbb{R}$ if $i \in I$.
\end{definition}

The (2-)Minkowski property was first considered in \cite{JO19}.
It follows from \Cref{def:mink} that every subgraph of a $d$-Minkowski graph is also $d$-Minkowski.
In fact, the $d$-Minkowski property is stronger than it first appears: any equivalent framework $(G,q)$ satisfying these properties is also \emph{non-congruent} to $(G,p)$.

\begin{lemma}\label{lem:minkisom}
    Let $p=\{p(v_1),\ldots,p(v_{d+1})\}$ be a set of generic real points in $\mathbb{R}^d$,
    and let $q=\{q(v_1),\ldots, q(v_{d+1})\}$ be a set of points in $\mathbb{R}^{d_1} \times (\mathfrak{i} \mathbb{R})^{d_2}$ for $d_1 \geq 0$ and $d_2 \geq 1$ with $d_1+d_2=d$.
    Then there is no pair $A \in O(d,\mathbb{C})$ and $x \in \mathbb{C}^d$ such that $p(v_i) = A q(v_i) +x$ for all $i \in [d+1]$.
\end{lemma}

\begin{proof}
    By translating the sets $p$ and $q$,
    we may suppose that $p(v_{d+1})=q(v_{d+1}) = \mathbf{0}$.
    Suppose for contradiction that there does exist some $A \in O(d,\mathbb{C})$ such that $p(v_i) = A q(v_i)$ for all $i \in [d]$.
    For each $i \in [d]$, fix the real vector $\tilde{q}(v_i) \in \mathbb{R}^d$ defined as follows:
    \begin{equation*}
        \tilde{q}_j(v_i) =
        \begin{cases}
            q_j(v_i) &\text{if } j \leq d_1, \\
            - \mathfrak{i} q_j(v_i) &\text{if } j \geq d_1 +1.
        \end{cases}
    \end{equation*}
    With this, we define the real $d \times d$ matrices $P = [p(v_1) \ \cdots \ p(v_d)]$ and $Q = [\tilde{q}(v_1) \ \cdots \ \tilde{q}(v_d)]$.
    Since $p$ is generic and $q$ is formed from $p$ by an invertible linear transformation, both matrices $P$ and $Q$ are invertible.
    Let $D$ be the $d \times d$ diagonal matrix where $D_{jj} = 1$ if $j \leq d_1$ and $D_{jj} = -1$ if $j \geq d_1 +1$.
    For each pair $i,j \in [d]$, we have that
    \begin{equation*}
        \tilde{q}(v_i)^T D \tilde{q}(v_j) = q(v_i)^T q(v_j) = (Ap(v_i))^T  (Ap(v_j)) = p(v_i)^T A^T A p(v_j) = p(v_i)^T p(v_j),
    \end{equation*}
    and so $P^T P = Q^T D Q$.
    However, Sylvester's law of inertia informs us that the eigenvalues of $P^T P$ are all positive, while $Q^T D Q$ has $d_1$ positive and $d_2 \geq 1$ negative eigenvalues, a contradiction.
\end{proof}

A $d$-Minkowski graph necessarily has different real and complex realisation number.

\begin{lemma}\label{lem:grnotmink}
    Let $G$ be a graph with at least $d+1$ vertices.
    If $r_d(G) = c_d(G)$, then $G$ is not $d$-Minkowski.
    In particular, if $G$ is globally $d$-rigid, then $G$ is not $d$-Minkowski.
\end{lemma}

\begin{proof}
    If $G$ is $d$-Minkowski, then \Cref{lem:minkisom} implies for any generic real $(G,p)$, there exists an equivalent non-congruent $(G,q)$ that is not real, and hence $r_d(G) \neq c_d(G)$.
    If $G$ is globally $d$-rigid,
    then $r_d(G)=c_d(G)$ by \cite[Thm 1]{Gortler2014},
    which hence implies $G$ is not $d$-Minkowski.    
\end{proof}

Our main result in this subsection is that graphs with odd $d$-realisation number must be $d$-Minkowski.

\begin{proposition}\label{prop:mink1}
    If $G$ is a $d$-rigid graph and $c_d(G)$ is odd but not equal to one,
    then $G$ is $d$-Minkowski.
\end{proposition}

\begin{proof}
    We first observe that $G$ must have at least $d+2$ vertices, since any $d$-rigid graph with $d+1$ or less vertices is a simplex and hence globally $d$-rigid.
    Choose any real generic $d$-dimensional realisation $p$ of $G$.
    Let
    \[
    S = \left\{ q \in (\mathbb{C}^d)^V \; : \; \|p (v)-p(w)\|^2 = \|q (v)-q(w)\|^2 \text{ for all } vw \in E\right\} \, .
    \]
    For any realisation $q \in S$,
    we define $q^*$ to be the realisation where we apply the complex conjugate to each coordinate of each vertex.
    For any edge $vw \in E$,
    we see that
    \begin{equation*}
        \|q^*(v) - q^*(w)\|^2 = \overline{\|q(v) - q(w)\|^2} = \overline{\|p(v) - p(w)\|^2} =\|p(v) - p(w)\|^2,
    \end{equation*}
    hence the map $q \mapsto q^*$ is an involution on $S$.
    As complex conjugation is an automorphism of $\CC$, it is straightforward to verify that $q_1 \sim q_2$ if and only if $q_1^* \sim q_2^*$.
    Explicitly, if $q_1 = Aq_2 + x$ for some $A \in O(d,\CC)$ and $x \in \CC^d$, then by automorphism properties we have $q_1^* = A^* q_2^* + x^*$ where $A^*$ is the entrywise complex conjugate of $A$.
    Furthermore, we have $A^* \in O(d,\CC)$ as $I_n = (I_n)^* = (A^\top A)^* = (A^*)^\top A^*$.
    Hence there exists a well-defined involution
    \begin{equation*}
        {}^* : C_d(G,p) \longrightarrow C_d(G,p), ~ [q] \longmapsto [q^*].
    \end{equation*}
    We now separate $C_d(G,p)$ into three families:
    \begin{align*}
        S_\mathbb{R} &:= \left\{ [q] \in C_d(G,p) : [q] \text{ contains a real realisation}  \right\}, \\
        S_M &:= \left\{ [q] \in C_d(G,p) : [q] \text{ does not contain a real realisation and } [q]^* = [q] \right\}, \\
        S_{\mathbb{C}} &:= C_d(G,p) \setminus \left( S_{\mathbb{R}} \cup S_M  \right).
    \end{align*}
    We note that $[q]^* = [q]$ for all $S_{\mathbb{R}}$ also.
    The involution maps every point in $S_\CC$ to another distinct point in $S_\CC$, implying that $|S_\CC|$ is even.
    As $G$ is not globally $d$-rigid and $p$ is generic, it follows from \cite[Lem.~2.38]{gortler2010characterizing} that $|S_{\mathbb{R}}|$ is even.
    Hence $|C_d(G,p)|$ is odd if and only if $|S_M|$ is odd.
    In particular, the set $S_M$ must be non-empty.
    
    Now choose $[q] \in C_d(G,p)$ such that $[q]^* = [q]$.
    Pick a representative $q$ of this equivalence class which is also contained in $X_{G,d}$ (which is guaranteed by \Cref{lem:quasixgd}),
    and fix the unique pair $A \in O(d, \mathbb{C})$, $x \in \mathbb{C}^d$ such that $q^*(v) = Aq(v) + x$ for each $v \in V$.
    As $q \in X_{G,d}$,
    it is easy to see that $x = 0$ and $A$ is a diagonal matrix (see, for example, the proof of \cite[Lem.~3.2]{dewar23}).
    As $A^T A = I_d$,
    there exists a set $I \subset [d]$ such that $A_{ii} = -1$ if $i \in I$ and $A_{ii} = 1$ otherwise.
    If $I = \emptyset$ then $[q] \in S_{\mathbb{R}}$.
    Thus if $[q] \in S_{M}$ then $I \neq \emptyset$,
    and our representative satisfies $q_i(v) \in \mathbb{R}$ if $i \in [d] \setminus I$ and $q_i(v) \in \mathfrak{i}\mathbb{R}$ if $i \in I$.
    This now concludes the proof.
\end{proof}

\subsection{Minimally isostatic graphs}

We recall that a $d$-rigid graph is a \emph{minimally $d$-isostatic graph} ($d$-MIG) if its only $d$-rigid subgraphs are complete on at most $d$ vertices.
As an example, $K_{3,3}$ is a 2-MIG as it is $2$-rigid and its only $2$-rigid subgraphs are single edges.
However, $K_4$ minus an edge is not a 2-MIG as it contains a copy of $K_3$.

A key step in the proof of \Cref{thm:mig} will be the following divisibility result:

\begin{theorem}[{\cite[Thm 1.1]{dewarsmithnixon2026}}]\label{thm:gcd}
    Let $G$ be a $d$-rigid graph on at least $d+1$ vertices.
    If $H$ is a spanning $d$-rigid subgraph of $G$,
    then $c_d(G) \mid c_d(H)$.
\end{theorem}

We are now ready to prove \Cref{thm:mig}, which we recall is as follows:

\mig*

\begin{proof}
    We prove this result by showing the contrapositive holds.
    Let $G$ be a minimally $d$-rigid graph that contains a proper minimally $d$-rigid subgraph $H$ with at least $d+1$ vertices.
    Fix $G' = G + K_{V(H)}$ where $K_{V(H)}$ is the complete graph on $V(H)$.
    
    We first show that $G'$ is not globally $d$-rigid, utilising the closure operator $\cl_d$ of the $d$-dimensional rigidity matroid $\cR_d$.
    Choose some edge $e \in E(G) \setminus E(H)$.
    Since $G$ is minimally $d$-rigid, $G -e$ is not $d$-rigid and hence $e \notin \cl_d(E(G) \setminus \{e\})$.
    As $H$ is $d$-rigid, we have
    \begin{equation*}
        \cl_d \left( E(G) \setminus \{e\} \right) = \cl_d \left( (E(G) \setminus \{e\}) \cup K_{V(H)}\right) = \cl_d \left( E(G') \setminus \{e\} \right),
    \end{equation*}
    hence $e$ is not contained in the closure of $E(G') \setminus \{e\}$.
    This implies $G' - e$ is not $d$-rigid, and therefore $G'$ is not redundantly $d$-rigid.
    By Hendrickson's theorem \cite{Hen92}, this implies $G'$ is not globally $d$-rigid.
    
    As $K_{V(H)}$ is a globally $d$-rigid graph on at least $d+1$ vertices, \Cref{lem:grnotmink} implies it is not $d$-Minkowski.
    As $G'$ contains $K_{V(H)}$, it follows that $G'$ is also not $d$-Minkowski.
    As $c_d(G') \neq 1$, \Cref{prop:mink1} gives that $c_d(G')$ is even.
    As $G$ is a $d$-rigid spanning subgraph of $G'$, it follows from \Cref{thm:gcd} that $c_d(G')|c_d(G)$, and hence $c_d(G)$ is also even.
\end{proof}

When $d=2$, we can drop the condition that $G$ is minimally $2$-rigid from \Cref{thm:mig}.
Jackson and Owen showed that the 2-realisation number of an $\cR_2$-connected graph $G$ is always a power of two and $c_2(G) = r_2(G)$ \cite[Thm 7.4]{JO19}.
Coupled with \Cref{lem:grnotmink}, this gives the following immediate corollary.

\begin{corollary}\label{cor:minkmconn}
    Let $G$ be an $\mathcal{R}_2$-connected graph.
    Then $G$ is not 2-Minkowski.
\end{corollary}

This allows us to deduce the following strengthening of \Cref{thm:mig} when $d=2$.

\begin{corollary}\label{cor:mig2d}
    Let $G$ be a 2-rigid graph with an odd number of realisations that is not 1. 
    Then $G$ is a 2-MIG.
\end{corollary}

\begin{proof}
    By \Cref{prop:mink1}, we have $G$ is 2-Minkowski.
    By \Cref{cor:minkmconn}, every circuit of $\cR_2$ is not 2-Minkowski.
    As the 2-Minkowski property is inherited by subgraphs, it follows that $G$ contains no circuits of $\cR_2$, and hence is minimally 2-rigid.
    The result now follows from \Cref{thm:mig}.
\end{proof}

Note that the converse is false; indeed, the bipartite graph $K_{3,3}$ is a 2-MIG but has exactly eight realisations (see, for example, \cite{capcodata}).

\begin{remark}\label{rem:globalrigid}
    It follows from \Cref{cor:mig2d} that if $G$ has an odd 2-realisation number then $G+e$ is globally 2-rigid for any choice of non-edge $e$.
    To see this, observe that every 2-MIG with at least 4 vertices is 3-connected.
    A simple edge-counting argument shows that if $G$ decomposes into $G_1 \cup G_2$ with $|V(G_1) \cap V(G_2)|=2$, 
    then minimal 2-rigidity implies that $|E(G_i)| = 2|V(G_i)|-3$ for some $i \in \{1,2\}$, forcing the subgraph $G_i$ to be minimally 2-rigid.
    Furthermore, the addition of any non-edge to a 2-MIG produces a 2-circuit, a consequence of \Cref{thm:2rigmatroid}.
    Applying \Cref{thm:2globrig} gives the result.
\end{remark}

We end this section by using the techniques developed to show the only 2-rigid planar graphs with odd 2-realisation number are globally 2-rigid.

\begin{proposition}\label{prop:planarodd}
    Let $G$ be a $2$-rigid planar graph.
    Then $c_2(G)$ is odd if and only if $G$ is globally 2-rigid.
\end{proposition}

\begin{proof}
    One direction is trivial, so assume $c_2(G)$ is odd and greater than 1.
    By \Cref{cor:mig2d},
    it suffices to assume that $G$ is minimally 2-rigid.
    Fix a planar embedding of $G$ to describe a face set $F$, and partition $F$ into sets $F_3, F_4, \ldots$ where $F_n$ is the set of all faces with $n$ sides.
    
    Suppose for contradiction that $G$ has no facial triangles (i.e.~$F_3=\emptyset$).
    By summing face-degrees, we see that
    \begin{equation}\label{eq:planarodd}
        2|E| = \sum_{n \geq 4} n|F_n| \geq 4\sum_{n \geq 4}|F_n| = 4|F|.
    \end{equation}
    By combining \eqref{eq:planarodd} with Euler's formula,
    we see that
    \begin{equation*}
        2 = |V| - |E| + |F| \leq |V| - |E| + \frac{2}{4}|E| = |V| - \frac{1}{2} |E| \quad \implies \quad |E| \leq 2|V| - 4 < 2|V|-3,
    \end{equation*}
    contradicting that $G$ is minimally 2-rigid.
    Hence $G$ contains a copy of $K_3$ and so is not a 2-MIG.
    The result now follows from \Cref{thm:mig}.
\end{proof}

We cannot improve \Cref{prop:planarodd} to higher genus graphs.
To see this, take the graph given in \Cref{fig:smallestodd} which has 2-realisation number of 45.
By \Cref{prop:planarodd},
our chosen graph is not planar.
As any 8-vertex graph with at most 21 edges can be embedded in the torus (see \cite[Thm 1]{Duke1972}),
the graph is toroidal.

\begin{remark}
    The family of planar 2-rigid graphs, here denoted $\mathcal{F}_{pl}$, has an additional property of significance:
    the integer set
    \begin{equation*}
        c_2( \mathcal{F}_{pl}) = \{ c_2(G) : G \in \mathcal{F}_{pl} \} \subset \mathbb{N}
    \end{equation*}
    is closed under multiplication.
    This follows from three key facts:
    (i) if we glue two 2-rigid graphs $G_1,G_2$ at a shared triangle (i.e.~$G_1 \cap G_2 \cong K_3$),
    then $c_2(G_1 \cup G_2) = c_2(G_1) c_2(G_2)$: this follows straightforwardly from the gluing lemma \cite[Lem.~11.1.9]{Whiteley96};
    (ii) every planar 2-rigid graph contains at least one triangle: this is shown in the proof of \Cref{prop:planarodd};
    and (iii) gluing two planar graphs at a facial triangle produces another planar graph.
    This argument similarly holds for any family of $d$-rigid graphs $\mathcal{F}$ where every member contains a copy of $K_{d+1}$ and gluing two graphs in $\mathcal{F}$ along some copy of $K_{d+1}$ produces another graph in $\mathcal{F}$.
\end{remark}

\section{An infinite family of graphs with odd realisation number} \label{sec:infiniteoddgraphs}

In this section we construct an infinite family of minimally 2-rigid graphs with odd 2-realisation number.
To do so, we need to introduce the following terminology from \cite{ghl24}.

\begin{definition}
    Let $G=(V,E)$ be a graph with edge $u_1u_2$ and vertex $x$ that is distinct from $u_1,u_2$.
    Then the triple $(G,u_1 u_2,x)$ is a \emph{calligraph} if:
    \begin{enumerate}
        \item $|E|=2|V|-4$, and,
        \item either $G+ x u_1$ or $G + x u_2$ is 2-rigid.
    \end{enumerate}
\end{definition}

\begin{definition}
    Let $G_1=(V_1,E_1)$ and $G_2=(V_2,E_2)$ be graphs with $V_1 \cap V_2 = \{x,u_1,u_2\}$ and $E_1 \cap E_2 = \{u_1 u_2\}$.
    Given $G = G_1 \cup G_2$,
    we say $(G_1,G_2,u_1 u_2, x)$ is a \emph{calligraphic split} of $G$ if $(G_i,u_1 u_2,x)$ is a calligraph for each $i \in \{1,2\}$.
\end{definition}

\begin{figure}[tp]
	\begin{center}
            \begin{tikzpicture}
            \begin{scope}[scale=0.75]
			\node[vertex] (u1) at (0,0) {};
			\node[vertex] (u2) at (2,0) {};
			\node[vertex] (x) at (1,1) {};
            
			\node at (-0.3,-0.3) {$u_1$};
			\node at (2.3,-0.3) {$u_2$};
			\node at (1,1.3) {$x$};
            
			\node at (1, -1) {$L$};

                \draw[edge] (u1)edge(u2);
			\draw[edge] (x)edge(u1);
            \end{scope}
            
            \begin{scope}[scale=0.75, xshift = 130]
			\node[vertex] (u1) at (0,0) {};
			\node[vertex] (u2) at (2,0) {};
			\node[vertex] (x) at (1,1) {};
            
			\node at (-0.3,-0.3) {$u_1$};
			\node at (2.3,-0.3) {$u_2$};
			\node at (1,1.3) {$x$};
            
			\node at (1, -1) {$R$};

                \draw[edge] (u1)edge(u2);
			\draw[edge] (x)edge(u2);
            \end{scope}
            
            \begin{scope}[scale=0.75, xshift = 260]
			\node[vertex] (u1) at (0,0) {};
			\node[vertex] (u2) at (2,0) {};
			\node[vertex] (v) at (1,1) {};
			\node[vertex] (x) at (1,2) {};
            
			\node at (-0.3,-0.3) {$u_1$};
			\node at (2.3,-0.3) {$u_2$};
			\node at (1.3,1.3) {$v$};
			\node at (1,2.3) {$x$};
            
			\node at (1, -1) {$C_v$};

                \draw[edge] (u1)edge(u2);
			\draw[edge] (x)edge(v);
			\draw[edge] (u1)edge(v);
			\draw[edge] (u2)edge(v);
            \end{scope}
    \end{tikzpicture}
	\end{center}
	\caption{The three base calligraphs used in \Cref{thm:calligraph}.}
    \label{fig:basecalligraph}
\end{figure}
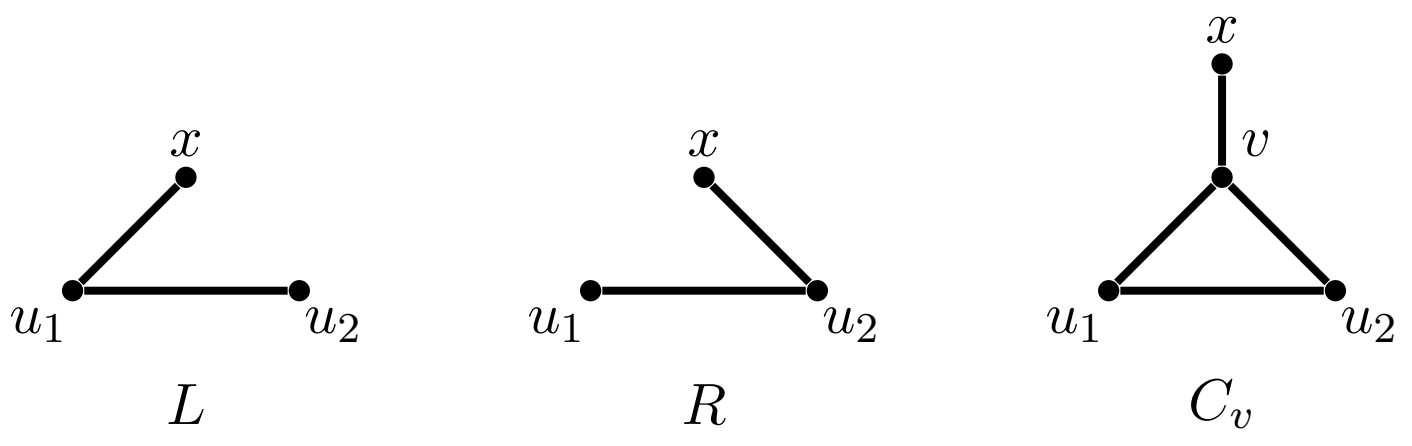

The key tool we will utilise is a method to inductively compute 2-realisation numbers via calligraphic splits.
This inductive method begins with the \emph{base calligraphs} $L, R$ and $C_v$ where $v \notin \{u_1,u_2,x\}$, defined in \Cref{fig:basecalligraph}.

\begin{theorem}[{\cite[Thm I]{ghl24}}]\label{thm:calligraph}
    Fix $\mathcal{C}$ to be the set of all graphs $G$ containing edge $u_1 u_2$ and vertex $x$ such that $(G,u_1u_2,x)$ is a calligraph.
    There exists a unique function $[ \cdot ] : \mathcal{C} \rightarrow \mathbb{Z}^3$ satisfying:
    \begin{enumerate}
        \item $[L] = (1,1,0)$, $[R] = (1,0,1)$ and $[C_v] = (2,0,0)$ for any choice of $v \notin \{u_1,u_2,x\}$;
        \item if $(G_1,G_2,u_1u_2,x)$ is a calligraphic split of a minimally 2-rigid graph $G$,
        then $c_2(G) = [G_1] \otimes [G_2]$,
        where $\otimes$ is the bilinear form
        \begin{equation*}
            (x,y,z) \otimes (x',y',z') = x x' - y y' - z z'.
        \end{equation*}
    \end{enumerate}
\end{theorem}

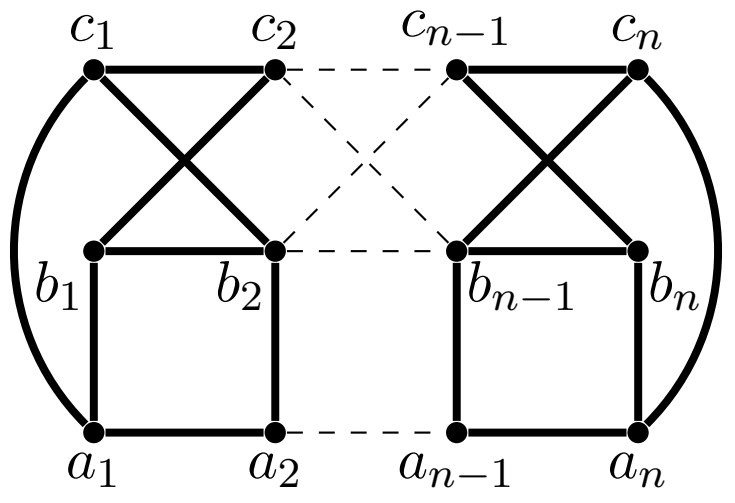
\begin{figure}
\centering
            \begin{tikzpicture}[scale=1.2]
			\node[vertex] (00) at (0,0) {};
            \node[below] at (0,0) {$a_1$};
			\node[vertex] (01) at (0,1) {};
            \node[below left=-2pt] at (0,1) {$b_1$};
			\node[vertex] (02) at (0,2) {};
            \node[above] at (0,2) {$c_1$};

			\node[vertex] (10) at (1,0) {};
            \node[below] at (1,0) {$a_2$};
			\node[vertex] (11) at (1,1) {};
            \node[below left=-2pt] at (1,1) {$b_2$};
			\node[vertex] (12) at (1,2) {};
            \node[above] at (1,2) {$c_2$};

			\node[vertex] (20) at (2,0) {};
            \node[below] at (2,0) {$a_{n-1}$};
			\node[vertex] (21) at (2,1) {};
            \node[below right=-2pt] at (2,1) {$b_{n-1}$};
			\node[vertex] (22) at (2,2) {};
            \node[above] at (2,2) {$c_{n-1}$};

            \node[vertex] (30) at (3,0) {};
            \node[below] at (3,0) {$a_{n}$};
			\node[vertex] (31) at (3,1) {};
            \node[below right=-2pt] at (3,1) {$b_n$};
			\node[vertex] (32) at (3,2) {};
            \node[above] at (3,2) {$c_n$};

			\draw[edge] (00)edge(10);
			\draw[black,dashed] (10)--(20);
			\draw[edge] (20)edge(30);

			\draw[edge] (01)edge(11);
			\draw[black,dashed] (11)--(21);
			\draw[edge] (21)edge(31);

			\draw[edge] (02)edge(12);
			\draw[black,dashed] (12)--(22);
			\draw[edge] (22)edge(32);

			\draw[edge] (00)edge(01);
			\draw[edge] (10)edge(11);
			\draw[edge] (20)edge(21);
			\draw[edge] (30)edge(31);

			\draw[edge] (01)edge(12);
			\draw[edge] (11)edge(02);
            
			\draw[black,dashed] (11)--(22);
			\draw[black,dashed] (21)--(12);

                \draw[edge] (21)edge(32);
			\draw[edge] (31)edge(22);

            \draw[edge] (00) to [bend left = 45] (02);

            \draw[edge] (30) to [bend right = 45] (32);
		\end{tikzpicture}
	\caption{A schematic of the family of minimally 2-rigid graphs $A_n$ for $n\geq 2$.}
    \label{fig:odd}
\end{figure}

We now define our key family of graphs: see \Cref{fig:odd} for an illustration.
\begin{definition}\label{def:An}
    For each $n \geq 2$, we define $A_n$ to be the graph with vertices and edges:
    \begin{align*}
    V(A_n) &= \{a_1,a_2,\dots, a_n,b_1,b_2,\dots, b_n, c_1,c_2,\dots, c_n\} \, ,\\
    E(A_n) &= \{a_ia_{i+1} \; : \; 1 \leq i \leq n-1\} \cup \{b_ib_{i+1} \; : \; 1 \leq i \leq n-1\} \cup \{c_ic_{i+1} \; : \; 1 \leq i \leq n-1\} \\
    &\; \cup \{a_ib_{i} \; : \; 1 \leq i \leq n\} \cup \{b_ic_{i+1} \; : \; 1 \leq i \leq n-1\} \cup \{c_ib_{i+1} \; : \; 1 \leq i \leq n-1\} \cup \{a_1c_1,a_nc_n\} \, .
    \end{align*} 
\end{definition}

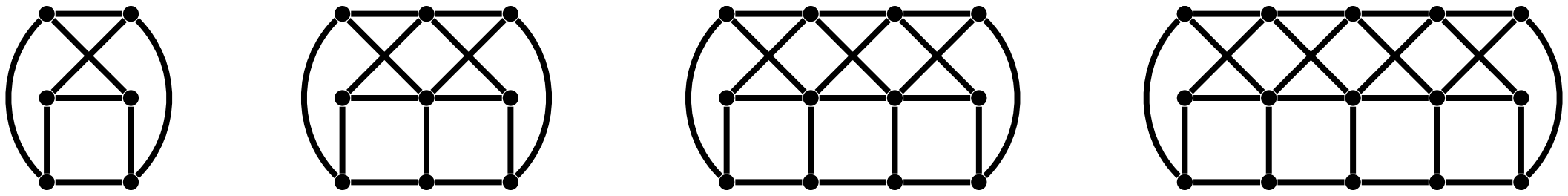
\begin{figure}[tp]
\centering
            \begin{tikzpicture}
            \begin{scope}[scale=0.75, xshift = -100]
			\node[vertex] (00) at (0,0) {};
			\node[vertex] (01) at (0,1) {};
			\node[vertex] (02) at (0,2) {};
                
			\node[vertex] (10) at (1,0) {};
			\node[vertex] (11) at (1,1) {};
			\node[vertex] (12) at (1,2) {};

			\draw[edge] (00)edge(10);

			\draw[edge] (01)edge(11);

			\draw[edge] (02)edge(12);
            
			\draw[edge] (00)edge(01);
			\draw[edge] (10)edge(11);

			\draw[edge] (01)edge(12);
			\draw[edge] (11)edge(02);

            \draw[edge] (00) to [bend left = 45] (02);

            \draw[edge] (10) to [bend right = 45] (12);
            \end{scope}

            \begin{scope}[scale=0.75]
			\node[vertex] (00) at (0,0) {};
			\node[vertex] (01) at (0,1) {};
			\node[vertex] (02) at (0,2) {};
                
			\node[vertex] (10) at (1,0) {};
			\node[vertex] (11) at (1,1) {};
			\node[vertex] (12) at (1,2) {};

			\node[vertex] (20) at (2,0) {};
			\node[vertex] (21) at (2,1) {};
			\node[vertex] (22) at (2,2) {};

			\draw[edge] (00)edge(10);
			\draw[edge] (10)edge(20);

			\draw[edge] (01)edge(11);
			\draw[edge] (11)edge(21);

			\draw[edge] (02)edge(12);
			\draw[edge] (12)edge(22);

			\draw[edge] (00)edge(01);
			\draw[edge] (10)edge(11);
			\draw[edge] (20)edge(21);

			\draw[edge] (01)edge(12);
			\draw[edge] (11)edge(02);
            
			\draw[edge] (11)edge(22);
			\draw[edge] (21)edge(12);

            \draw[edge] (00) to [bend left = 45] (02);

            \draw[edge] (20) to [bend right = 45] (22);
            \end{scope}
            
            \begin{scope}[scale=0.75, xshift = 130]
			\node[vertex] (00) at (0,0) {};
			\node[vertex] (01) at (0,1) {};
			\node[vertex] (02) at (0,2) {};
			\node[vertex] (02) at (0,2) {};

			\node[vertex] (10) at (1,0) {};
			\node[vertex] (11) at (1,1) {};
			\node[vertex] (12) at (1,2) {};

			\node[vertex] (20) at (2,0) {};
			\node[vertex] (21) at (2,1) {};
			\node[vertex] (22) at (2,2) {};

                \node[vertex] (30) at (3,0) {};
			\node[vertex] (31) at (3,1) {};
			\node[vertex] (32) at (3,2) {};

			\draw[edge] (00)edge(10);
			\draw[edge] (10)edge(20);
			\draw[edge] (20)edge(30);

			\draw[edge] (01)edge(11);
			\draw[edge] (11)edge(21);
			\draw[edge] (21)edge(31);

			\draw[edge] (02)edge(12);
			\draw[edge] (12)edge(22);
			\draw[edge] (22)edge(32);

			\draw[edge] (00)edge(01);
			\draw[edge] (10)edge(11);
			\draw[edge] (20)edge(21);
			\draw[edge] (30)edge(31);

			\draw[edge] (01)edge(12);
			\draw[edge] (11)edge(02);
            
			\draw[edge] (11)edge(22);
			\draw[edge] (21)edge(12);

                \draw[edge] (21)edge(32);
			\draw[edge] (31)edge(22);

            \draw[edge] (00) to [bend left = 45] (02);

            \draw[edge] (30) to [bend right = 45] (32);
            \end{scope}

            \begin{scope}[scale=0.75, xshift = 285]
			\node[vertex] (00) at (0,0) {};
			\node[vertex] (01) at (0,1) {};
			\node[vertex] (02) at (0,2) {};
			\node[vertex] (02) at (0,2) {};

			\node[vertex] (10) at (1,0) {};
			\node[vertex] (11) at (1,1) {};
			\node[vertex] (12) at (1,2) {};

			\node[vertex] (20) at (2,0) {};
			\node[vertex] (21) at (2,1) {};
			\node[vertex] (22) at (2,2) {};

                \node[vertex] (30) at (3,0) {};
			\node[vertex] (31) at (3,1) {};
			\node[vertex] (32) at (3,2) {};

                \node[vertex] (40) at (4,0) {};
			\node[vertex] (41) at (4,1) {};
			\node[vertex] (42) at (4,2) {};

			\draw[edge] (00)edge(10);
			\draw[edge] (10)edge(20);
			\draw[edge] (20)edge(30);
			\draw[edge] (30)edge(40);

			\draw[edge] (01)edge(11);
			\draw[edge] (11)edge(21);
			\draw[edge] (21)edge(31);
			\draw[edge] (31)edge(41);

			\draw[edge] (02)edge(12);
			\draw[edge] (12)edge(22);
			\draw[edge] (22)edge(32);
			\draw[edge] (32)edge(42);

			\draw[edge] (00)edge(01);
			\draw[edge] (10)edge(11);
			\draw[edge] (20)edge(21);
			\draw[edge] (30)edge(31);
			\draw[edge] (40)edge(41);

			\draw[edge] (01)edge(12);
			\draw[edge] (11)edge(02);
            
			\draw[edge] (11)edge(22);
			\draw[edge] (21)edge(12);

                \draw[edge] (21)edge(32);
			\draw[edge] (31)edge(22);

                \draw[edge] (31)edge(42);
			\draw[edge] (41)edge(32);

            \draw[edge] (00) to [bend left = 45] (02);

            \draw[edge] (40) to [bend right = 45] (42);
            \end{scope}
            \end{tikzpicture}
        \caption{The first four graphs in the family $A_n$. For $n \geq3$, these have odd 2-realisation numbers. From left to right: $c_2(A_2) = 8$, $c_2(A_3) = 103$, $c_2(A_4) = 1111$, $c_2(A_5) = 11159$.}
        \label{fig:odd_25}
\end{figure}

The first observation about the family $A_n$ is that they are all minimally 2-rigid.
\begin{lemma} \label{lem:An-rigid}
    For each $n \geq 2$, the graph $A_n$ is minimally 2-rigid.
\end{lemma}
\begin{proof}
Note that $A_2$ is isomorphic to $K_{3,3}$, which is well known to be minimally 2-rigid. For all $n\geq 3$, $A_n$ can be obtained from $A_{n-1}$ by the following sequence of extension moves:
\begin{align*}
    \text{1-extension:}& \quad \text{add $c_n$ adjacent to $a_{n-1},b_{n-1},c_{n-1}$, remove $a_{n-1}c_{n-1}$,}\\
    \text{0-extension:}& \quad \text{add $b_n$ adjacent to $b_{n-1},c_{n-1}$,}\\
    \text{1-extension:}& \quad \text{add $a_n$ adjacent to $a_{n-1},b_{n},c_{n}$, remove $a_{n-1}c_{n}$}.
\end{align*}
As these extension moves preserve minimal $2$-rigidity by \Cref{lem:01ext}, each graph $A_n$ is minimally 2-rigid.
\end{proof}

Consider the first four graphs in the family $A_n$, depicted in \Cref{fig:odd_25}.
The 2-realisation numbers of the first three graphs are given in the database \cite{capcodata}.
The 2-realisation number of the fourth can be computed via the algorithm in \cite{cggkls}.
Barring $A_2$, each graph has an odd 2-realisation number.
Our main result of this section is that this pattern continues: that $c_2(A_n)$ is odd for all $n \geq 3$.

From now on, we fix $u_1 = a_1$, $u_2 = b_1$ and $x = c_1$ whenever we use calligraphs for the remainder of the section.
We also let $\widetilde{A}_n := A_n - a_1c_1$.
With this, each triple $(\widetilde{A}_n, a_1 b_1, c_1)$ is a calligraph.
We define the sequence of vectors $((x_n,y_n,z_n))_{n \geq 2}$ in $\mathbb{Z}^3$ such that $[\widetilde{A}_n] = (x_n,y_n,z_n)$ for each $n \geq 2$.

\begin{lemma}\label{lem:smallvecs}
    $[\widetilde{A}_2] = (11,3,3)$ and $[\widetilde{A}_3] = (122,19,58)$.
\end{lemma}

\begin{proof}
    Throughout the proof, we will give the 2-realisation numbers of some specific minimally 2-rigid graphs without proof: these can be computed using the algorithm in \cite{cggkls} or looked up in the database given in \cite{capcodata}.
    
    First, set $n = 2$.
    \Cref{thm:calligraph} gives us the equations
    \begin{align*}
        8 &= c_2(\widetilde{A}_2 \cup L) = [\widetilde{A}_2] \otimes [L] = x_2 - y_2, \\
        8 &= c_2(\widetilde{A}_2 \cup R) = [\widetilde{A}_2] \otimes [R] = x_2 - z_2, \\
        22 &= c_2(\widetilde{A}_2 \cup C_v) = [\widetilde{A}_2] \otimes [C_v] = 2 x_2,
    \end{align*}
    which have the unique solution $(11,3,3)$.

    Now set $n = 3$.
    By computing the realisation numbers directly,
    we obtain the equations
    \begin{align*}
        103 &= c_2(\widetilde{A}_3 \cup L) = [\widetilde{A}_3] \otimes [L] = x_3 - y_3, \\
        64 &= c_2(\widetilde{A}_3 \cup R) = [\widetilde{A}_3] \otimes [R] = x_3 - z_3, \\
        244 &= c_2(\widetilde{A}_3 \cup C_v) = [\widetilde{A}_3] \otimes [C_v] = 2 x_3,
    \end{align*}
    which have the unique solution $(122,19,58)$.
\end{proof}

\begin{figure}[tp]
	\begin{center}
            \begin{tikzpicture}
            \begin{scope}[scale=0.75]
			\node[vertex] (00) at (0,0) {};
			\node[vertex] (01) at (0,1) {};
			\node[vertex] (02) at (0,2) {};

			\node[vertex] (10) at (1,0) {};
			\node[vertex] (11) at (1,1) {};
			\node[vertex] (12) at (1,2) {};

			\node[vertex] (20) at (2,0) {};
			\node[vertex] (21) at (2,1) {};
			\node[vertex] (22) at (2,2) {};

                \node[vertex] (30) at (3,0) {};
			\node[vertex] (31) at (3,1) {};
			\node[vertex] (32) at (3,2) {};

			\draw[edge] (00)edge(10);
			\draw[black,dashed] (10)--(20);
			\draw[edge] (20)edge(30);

			\draw[edge] (01)edge(11);
			\draw[black,dashed] (11)--(21);
			\draw[edge] (21)edge(31);

			\draw[edge] (02)edge(12);
			\draw[black,dashed] (12)--(22);
			\draw[edge] (22)edge(32);

			\draw[edge] (00)edge(01);
			\draw[edge] (10)edge(11);
			\draw[edge] (20)edge(21);
			\draw[edge] (30)edge(31);

			\draw[edge] (01)edge(12);
			\draw[edge] (11)edge(02);
            
			\draw[black,dashed] (11)--(22);
			\draw[black,dashed] (21)--(12);

                \draw[edge] (21)edge(32);
			\draw[edge] (31)edge(22);

            \draw[edge] (01)edge(02);

            \draw[edge] (30) to [bend right = 45] (32);

            \node at (1.5,-1) {$\widetilde{A}_n \cup R$};
            \end{scope} 
            
            \begin{scope}[scale=0.75, xshift=200]
			\node[vertex] (00) at (0,0) {};
			\node[vertex] (01) at (0,1) {};
			\node[vertex] (02) at (0,2) {};

			\node[vertex] (10) at (1,0) {};
			\node[vertex] (11) at (1,1) {};
			\node[vertex] (12) at (1,2) {};

			\node[vertex] (20) at (2,0) {};
			\node[vertex] (21) at (2,1) {};
			\node[vertex] (22) at (2,2) {};

                \node[vertex] (30) at (3,0) {};
			\node[vertex] (31) at (3,1) {};
			\node[vertex] (32) at (3,2) {};
            
			\node[vertex] (v) at (-1,1) {};

			\draw[edge] (00)edge(10);
			\draw[black,dashed] (10)--(20);
			\draw[edge] (20)edge(30);

			\draw[edge] (01)edge(11);
			\draw[black,dashed] (11)--(21);
			\draw[edge] (21)edge(31);

			\draw[edge] (02)edge(12);
			\draw[black,dashed] (12)--(22);
			\draw[edge] (22)edge(32);

			\draw[edge] (00)edge(01);
			\draw[edge] (10)edge(11);
			\draw[edge] (20)edge(21);
			\draw[edge] (30)edge(31);

			\draw[edge] (01)edge(12);
			\draw[edge] (11)edge(02);
            
			\draw[black,dashed] (11)--(22);
			\draw[black,dashed] (21)--(12);

                \draw[edge] (21)edge(32);
			\draw[edge] (31)edge(22);

            \draw[edge] (v)edge(00);
            \draw[edge] (v)edge(01);
            \draw[edge] (v)edge(02);

            \draw[edge] (30) to [bend right = 45] (32);

            \node at (1.5,-1) {$\widetilde{A}_n \cup C_v$};
            \end{scope} 
		\end{tikzpicture}
	\end{center}
	\caption{The graphs $\widetilde{A}_n \cup R = \widetilde{A}_n + b_1 c_1$ and $\widetilde{A}_n \cup C_v$.}
    \label{fig:RandCv}
\end{figure}
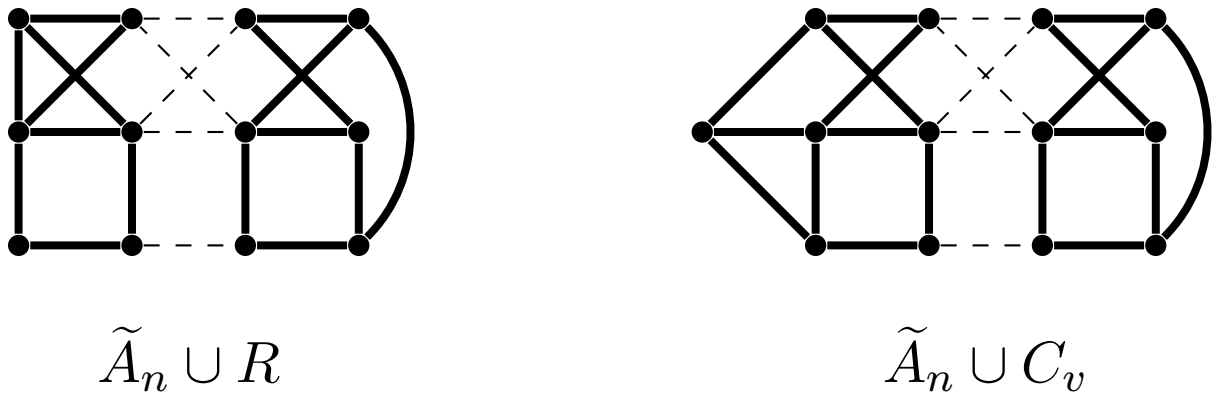

We now use these base cases to compute $[\widetilde{A}_n] = (x_n,y_n,z_n)$ inductively.
We will utilise the fact that we have the following calligraphic splits of $A_n$ for each $s,t \geq 2$ with $s+t-1=n$:
\begin{equation}\label{eq:addingtheodds}
    \widetilde{A}_n \cup L = A_n \qquad \text{and} \qquad \widetilde{A}_s \cup \widetilde{A}_t \cong A_n.
\end{equation}
Two other graphs of interest we will use are $\widetilde{A}_n \cup R = \widetilde{A}_n + b_1 c_1$ and $\widetilde{A}_n \cup C_v$, as pictured in \Cref{fig:RandCv}.
Both graphs are minimally 2-rigid;
the former can be constructed from a triangle by a sequence of 0-extensions, while the latter is formed from $A_n$ by a 1-extension.

\begin{lemma}\label{lem:allodd}
    For each $n \geq 3$,
    $x_n \equiv z_n \equiv 0 \mod 2$ and $y_n \equiv 1 \mod 2$.
\end{lemma}

\begin{proof}
By \Cref{thm:calligraph} and \eqref{eq:addingtheodds},
we have for every integer $n \geq 2$,
that
\begin{equation}\label{eq:Lmult}
    c_2(A_n) = [\widetilde{A}_n] \otimes [L] = x_n - y_n,
\end{equation}
and for every pair of integers $s,t \geq 2$ with $s+t -1 = n$,
we have
\begin{equation}\label{eq:genmult}
    c_2(A_n) = [\widetilde{A}_s] \otimes [\widetilde{A}_t] 
\end{equation}
Also, recall that $\widetilde{A}_n \cup R$ can be constructed via 2-dimensional 0-extensions from a triangle, hence it follows from \Cref{lem:0ext} that $c_2(\widetilde{A}_n \cup R) = 2^{3n-3}$.
\Cref{thm:calligraph} then implies that
\begin{equation}\label{eq:Rmult}
    c_2(\widetilde{A}_n \cup R) = [\widetilde{A}_n] \otimes [R] = 2^{3n-3}.
\end{equation}
Using these equations, we see that the following facts hold for all $n \geq 3$:
\begin{enumerate}
    \item \label{fact2} By \eqref{eq:Lmult},
    we have that $x_n - y_n = c_2(A_n)$.
    \item \label{fact3} By \eqref{eq:genmult} with $s = n-1$ and $t = 2$, and using $(x_2,y_2,z_2) = (11,3,3)$ from \Cref{lem:smallvecs},
    we have that $11x_{n-1} - 3y_{n-1} - 3z_{n-1} = c_2(A_n)$.
    \item \label{fact4} By \eqref{eq:genmult} with $s = n-2$ and $t = 3$, and using  $(x_3,y_3,z_3) = (122,19,58)$ from \Cref{lem:smallvecs},
    we have that $122x_{n-2} - 19 y_{n-2} - 58 z_{n-2} = c_2(A_n)$.
    \item \label{fact1} By \eqref{eq:Rmult},
    we have that $x_n - z_n = 2^{3n-3}$;
    hence, $x_n \equiv z_n \mod 2$ for all $n \geq 2$.
\end{enumerate}
We now proceed with an induction argument.
The statement holds when $n=3$ by \Cref{lem:smallvecs}.
Suppose now that $n \geq 4$ and for each integer $k$ satisfying $3 \leq k \leq n-1$ we have $x_k \equiv z_k \equiv 0 \mod 2$ and $y_k \equiv 1 \mod 2$.
By \ref{fact3},
we have that
\begin{equation*}
    c_2(A_n) \equiv 11 x_{n-1} - 3y_{n-1} - 3z_{n-1} \equiv 11 \cdot 0 - 3 \cdot 1 - 3 \cdot 0 \equiv 1 \mod 2.
\end{equation*}
It now follows from \ref{fact2} that $x_n - y_n \equiv 1 \mod 2$.
Combining this with \ref{fact1},
we either have that $x_n \equiv z_n \equiv 0 \mod 2$ and $y_n \equiv 1 \mod 2$ (which we want),
or $x_n \equiv z_n \equiv 1 \mod 2$ and $y_n \equiv 0 \mod 2$.
Using \ref{fact4} with $n+1$ instead of $n$,
we have that
\begin{equation*}
    c_2(A_{n+1}) \equiv 122 x_{n-1} - 19 y_{n-1} - 58 z_{n-1} \equiv 122 \cdot 0 - 19 \cdot 1 - 58 \cdot 0 \equiv 1 \mod 2.
\end{equation*}
Next,
using \ref{fact3} with $n+1$ instead of $n$,
we have that
\begin{equation*}
    1 \equiv c_2(A_{n+1}) \equiv 11 x_{n} - 3 y_{n} - 3 z_{n} \equiv (x_n - z_n) + y_n \equiv y_n \mod 2.
\end{equation*}
Hence we have the desired result that $x_n \equiv z_n \equiv 0 \mod 2$ and $y_n \equiv 1 \mod 2$.
\end{proof}

\infinite*

\begin{proof}
By \Cref{lem:allodd} and \eqref{eq:Lmult},
we have $c_2(A_n) \equiv x_n - y_n \equiv 1 \mod 2$ holds for all $n \geq 3$.
Hence, $c_2(A_n)$ is odd for all $n \geq 3$.
\end{proof}

\section{A non-constructive infinite family of graphs with odd realisation number}\label{sec:nonconstruct}

In this section, we introduce a gluing operation in \Cref{thm:oddglue} that takes as input two minimally 2-rigid graphs with odd realisation number satisfying certain conditions, and outputs two new minimally 2-rigid graphs, exactly one of which has odd realisation number.
This gluing operation is then used to give a non-constructive proof for \Cref{thm:allodd}.
As with the family of graphs introduced in \Cref{sec:infiniteoddgraphs}, this gluing operation relies on calligraphs and calligraphic splits introduced in \cite{ghl24}.

\begin{theorem} \label{thm:oddglue}
    Let $G_1 = (V_1, E_1)$ and $G_2 = (V_2,E_2)$ be graphs with $V_1 \cap V_2 = \{x,u_1,u_2\}$ and $E_1 \cap E_2 = \{u_1u_2\}$ such that $xu_1 \notin E_1 \cup E_2$.
    Suppose that:
    \begin{enumerate}
        \item $c_2(G_1 + xu_1)$ and $c_2(G_2 + xu_1)$ are odd and not equal to one;
        \item $c_2(G_1 + xu_2)$ and $c_2(G_2 + xu_2)$ are not odd. \label{cond:xu2-not-odd}
    \end{enumerate}
    Let $\widetilde{G}_2$ be the graph isomorphic to $G_2$ obtained by swapping the vertex labels of $u_1$ and $u_2$.
    Let $H = G_1 \cup G_2$ and $\widetilde{H} = G_1 \cup \widetilde{G}_2$.
    Then both $H$ and $\widetilde{H}$ are minimally 2-rigid, and exactly one of the integers $c_2(H)$ and $c_2(\widetilde{H})$ is odd and not equal to one.
\end{theorem}

\begin{proof}
    As both $c_2(G_1 + xu_1)$ and $c_2(G_2 + xu_1)$ are odd and not equal to one, it follows from \Cref{cor:mig2d} that both $G_1 + xu_1$ and $G_2 + xu_1$ are $2$-MIGs and hence minimally 2-rigid.
    As $xu_1 \notin E_1 \cup E_2$, it follows from \Cref{thm:2rigmatroid}\ref{thm:2rigmatroid:rigid} that $E_i = 2|V_i| - 4$ for each $i \in \{1,2\}$, hence both $G_1$ and $G_2$ are calligraphs.
    As the relabelling does not affect rigidity, it follows that $\widetilde{G}_2$ is also a calligraph.

    Note that Condition \ref{cond:xu2-not-odd} makes no assumption on $c_2(G_1 + xu_2)$ and $c_2(G_2 + xu_2)$ being finite, just that they are not odd.
    We now show that they are indeed finite, and in fact even.

    \begin{claim}
        $G_1 + xu_2$ and $G_2 + xu_2$ are minimally 2-rigid and $c_2(G_1 + xu_2)$ and $c_2(G_2 + xu_2)$ are even.
    \end{claim}

    \begin{claimproof}
        First note that $xu_2 \notin E_1$.
        Indeed, if $xu_2 \in E_1$ then $G_1 + xu_1$ would contain a triangle, contradicting that $G_1 + xu_1$ is a 2-MIG.
        The same reasoning gives that $xu_2 \notin E_2$.
        As $G_1+xu_1$ is a 2-MIG, it follows from \Cref{rem:globalrigid} that $G_1 + xu_1+ xu_2$ is globally 2-rigid.
        As globally 2-rigid graphs are redundantly rigid by \Cref{thm:2globrig}, it follows that $G_1 + xu_2$ is 2-rigid.
        Minimal 2-rigidity now follows from \Cref{thm:2rigmatroid}\ref{thm:2rigmatroid:rigid} and the fact that $G_1 + xu_2$ has $2|V_1|-3$ edges.
        The same reasoning gives that $G_2 +xu_2$ is also minimally 2-rigid.
        This implies that both $c_2(G_1 + xu_2)$ and $c_2(G_1 + xu_2)$ must be finite but not equal to one.
        Condition \ref{cond:xu2-not-odd} now implies they must both be even.    
        \end{claimproof}
    
    \begin{claim}
        $H$ and $\widetilde{H}$ are minimally 2-rigid.
    \end{claim}
    \begin{claimproof}
    We show that $H = (V,E)$ is minimally 2-rigid via \Cref{thm:2rigmatroid}\ref{thm:2rigmatroid:rigid}; the claim for $\widetilde{H}$ is identical.
    Recall that $E$ is independent in $\cR_2$ if and only if for any subgraph $H' = (V',E')$ of $H$, we have $|E'| \leq 2|V'| - 3$.
    It suffices to check this for induced subgraphs of $H$, so we will assume $H'$ is induced.
    Write $G_i' = (V_i',E_i')$ where $V_i' = V_i \cap V'$ and $E_i' = E_i \cap E'$ for each $i =1,2$.
    Note that each $E_i' \subseteq E_i$ is independent in $\cR_2$ by \Cref{thm:2rigmatroid} and that $G_i + xu_1$ is minimally 2-rigid.
    
    We proceed by case analysis on $V_1' \cap V_2'$.
    If $|V_1' \cap V_2'| \leq 1$, then $E'_1 \cap E'_2 = \emptyset$.
    Therefore
    \[
    |E'| = |E_1'| + |E_2'| \leq (2|V_1'| - 3) + (2|V_2'| - 3) = 2(|V'| + |V_1' \cap V_2'|) - 6 \leq 2|V'| - 4 \, .
    \]
    If $V_1' \cap V_2' = \{x,u_1,u_2\}$, then $E'_1 \cap E'_2 = \{u_1u_2\}$.
    As $G_i + xu_1$ is minimally 2-rigid, we have $E_i + xu_1$ is independent in $\cR_2$ and hence $|E'_i| \leq 2|V'_i| - 4$.
    Therefore
    \[
    |E'| = |E_1'| + |E_2'| - 1 \leq (2|V_1'| - 4) + (2|V_2'| - 4) -1 = 2(|V'| + |V_1' \cap V_2'|) - 9 \leq 2|V'| - 3 \, .
    \]
    If $V_1' \cap V_2' = \{x,u_j\}$ for some $j=1,2$, then $E'_1 \cap E'_2 = \emptyset$.
    Again, $G_i + xu_j$ being minimally 2-rigid implies that $|E'_i| \leq 2|V'_i| - 4$.
    Therefore
    \[
    |E'| = |E_1'| + |E_2'| \leq (2|V_1'| - 4) + (2|V_2'| - 4) = 2(|V'| + |V_1' \cap V_2'|) - 8 \leq 2|V'| - 4 \, .
    \]
    Finally, if $V_1' \cap V_2' = \{u_1,u_2\}$ then $E'_1 \cap E'_2 = \{u_1u_2\}$.
    In particular, $G'_1 \cap G'_2$ is a single edge and hence 2-rigid.
    It now follows from the gluing lemma \cite[Lemma 3.1.4]{Whiteley96} that $E'$ is independent in $\cR_2$.

    Finally, a similar counting argument gives that
    \[
    |E| =|E_1| + |E_2| - |E_1 \cap E_2| = (2|V_1| - 4) + (2|V_2| - 4) -1 = 2(|V| + |V_1 \cap V_2|) - 9 = 2|V| - 3 \, .
    \]
    Hence \Cref{thm:2rigmatroid}\ref{thm:2rigmatroid:rigid} implies $H$ is minimally 2-rigid.
    \end{claimproof}

    With this, we have the following calligraphic splits:
    \begin{enumerate}
        \item $(G_1,L,u_1u_2,x)$ is a calligraphic split of $G_1 + xu_1$; \label{split1}
        \item $(G_2,L,u_1u_2,x)$ is a calligraphic split of $G_2 + xu_1$, and \label{split2}\\ 
         $(\widetilde{G}_2,R,u_1u_2,x)$ is a calligraphic split of $\widetilde{G}_2 + xu_2 \cong G_2 + xu_1$;
        \item $(G_1,R,u_1u_2,x)$ is a calligraphic split of $G + xu_2$;\label{split3}
        \item $(G_2,R,u_1u_2,x)$ is a calligraphic split of $G_2 +xu_2$, and \label{split4}\\ 
        $(\widetilde{G}_2,L,u_1u_2,x)$ is a calligraphic split of $\widetilde{G}_2 + xu_1\cong G_2 +xu_2$.
        \item $(G_1, G_2, u_1u_2,x)$ is a calligraphic split of $H$, and \label{split5}\\ 
        $(G_1, \widetilde{G}_2, u_1u_2,x)$ is a calligraphic split of $\widetilde{H}$.
    \end{enumerate}
    Set $[G_1] = (a_1,b_1,c_1)$ and $[G_2]= (a_2,b_2,c_2)$.
    \begin{claim} \label{claim:vertex-swap}
        $[\widetilde{G}_2]= (a_2,c_2,b_2)$.
    \end{claim}
    \begin{claimproof}
        Let $[\widetilde{G}_2] = (\tilde{a}_2, \tilde{b}_2, \tilde{c}_2)$.
        Let $\Gamma$ be the graph obtained from $G_2$ by adding a new vertex $v$ and the three edges $xv,u_1v,u_2v$.
        Observe that $\Gamma$ is minimally 2-rigid, as it can be obtained from the minimally 2-rigid graph $G_2 + xu_1$ via a 1-extension.
        By the symmetry of the base graph $C_v$, both $(G_2,C_v,u_1u_2,x)$ and $(\widetilde{G}_2,C_v,u_1u_2,x)$ are calligraphic splits of $\Gamma$, hence applying \Cref{thm:calligraph} we obtain
        \[
        2a_2 = [G_2] \otimes [C_v] = [\widetilde{G}_2] \otimes [C_v] = 2 \tilde{a}_2 \, .
        \]
        
        Now consider the calligraphic splits \ref{split2} and \ref{split4}.
        As each case gives a calligraphic split of the same graph using $G_2$ and $\widetilde{G_2}$, we have
        \begin{align*}
            a_2 - b_2 = [G_2] \otimes [L] = [\widetilde{G}_2] \otimes [R] = \tilde{a}_2 - \tilde{c}_2 \, , \\
            a_2 - c_2 = [G_2] \otimes [R] = [\widetilde{G}_2] \otimes [L] = \tilde{a}_2 - \tilde{b}_2 \, .
        \end{align*}
        As $a_2 = \tilde{a}_2$, the claim follows.
    \end{claimproof}
    
    Applying \Cref{thm:calligraph} to our five calligraphic splits, we obtain the following:
    \begin{enumerate}
        \item $c_2(G_1+xu_1) = a_1-b_1$, and thus exactly one of $a_1$ and $b_1$ is odd;
        \item $c_2(G_2+xu_1) = a_2-b_2$, and thus exactly one of $a_2$ and $b_2$ is odd; 
        \item $c_2(G_1 + xu_2) = a_1-c_1$, and thus $a_1$ and $c_1$ have the same parity;
        \item $c_2(G_2 + xu_2) = a_2-c_2$, and thus $a_2$ and $c_2$ have the same parity;
        \item $c_2(H) = a_1a_2 -b_1b_2 - c_1c_2$ and $c_2(\widetilde{H}) = a_1a_2-b_1c_2 - c_1b_2$;
    \end{enumerate}
    where the final case utilises \Cref{claim:vertex-swap}.
    We now have four possible options to consider:
    \begin{itemize}
        \item $a_1 \equiv a_2 \equiv 0 \mod 2$:
        here we have that $b_1 \equiv b_2 \equiv 1 \mod 2$ and $c_1 \equiv c_2 \equiv 0 \mod 2$.
        Thus $c_2(H) \equiv 1 \mod 2$ and $c_2(\widetilde{H}) \equiv 0 \mod 2$.
        \item $a_1 \equiv a_2 \equiv 1 \mod 2$:
        here we have that $b_1 \equiv b_2 \equiv 0 \mod 2$ and $c_1 \equiv c_2 \equiv 1 \mod 2$.
        Thus $c_2(H) \equiv 0 \mod 2$ and $c_2(\widetilde{H}) \equiv 1 \mod 2$.
        \item $a_1 \equiv 0 \mod 2$ and $a_2 \equiv 1 \mod 2$:
        here we have that $b_1 \equiv c_2 \equiv 1 \mod 2$ and $c_1 \equiv b_2 \equiv 0 \mod 2$.
        Thus $c_2(H) \equiv 0 \mod 2$ and $c_2(\widetilde{H}) \equiv 1 \mod 2$.
        \item $a_1 \equiv 1 \mod 2$ and $a_2 \equiv 0 \mod 2$:
        here we have that $b_1 \equiv c_2 \equiv 0 \mod 2$ and $c_1 \equiv b_2 \equiv 1 \mod 2$.
        Thus $c_2(H) \equiv 0 \mod 2$ and $c_2(\widetilde{H}) \equiv 1 \mod 2$.
    \end{itemize}
    Each option shows that exactly one of $c_2(H)$ and $c_2(\widetilde{H})$ is odd.
    As both are minimally 2-rigid, they are not globally 2-rigid and hence neither have 2-realisation number equal to one.
\end{proof}

\Cref{thm:oddglue} gives rise to an alternative proof of \Cref{thm:allodd} as follows: 
\begin{itemize}
    \item Take any graph with at least two disjoint induced 4-cycles and an odd 2-realisation number; for example, the graph pictured in \Cref{fig:smallestodd}.
    \item Label one of the induced 4-cycles by $(x,u_1,u_2,y)$, and let $G_1$ be the graph with $xu_1$ removed.
    \item The graph $G_1 + xu_2$ has $2|V|-3$ edges and at least four vertices, and so cannot be globally 2-rigid by \Cref{thm:2globrig}. 
    Since it contains the triangle $\{x,y,u_2\}$, \Cref{cor:mig2d} implies that $G_1 + xu_2$ cannot have odd 2-realisation number.
    \item We now let $G_2 \cong G_1$ with $x$, $u_1$ and $u_2$ the only shared vertex labels, and apply \Cref{thm:oddglue}.
\end{itemize}
Since both of the graphs $H,\widetilde{H}$ contain at least two disjoint induced 4-cycles and exactly one of them has an odd realisation number greater than 1, we can keep repeating this process indefinitely to obtain an infinite family of minimally 2-rigid graphs with an odd 2-realisation number.
We emphasise that this procedure is non-constructive, as we do not know which of $H$ or $\widetilde{H}$ has odd 2-realisation number, just that one of them must do.

\begin{remark}
    We briefly note that the family of graphs $A_n$ from \Cref{def:An} can be obtained via the gluing procedure in \Cref{thm:oddglue}, but emphasise that it cannot alone deduce that $c_2(A_n)$ is odd.
    
    To see this, suppose we have the induction hypothesis that $c_2(A_s)$ is odd for all $3 \leq s < n$.
    In the notation of \Cref{thm:oddglue}, let $a_1 = u_1$, $b_1 = u_2$ and $c_1 = x$, and set $G_1 = \widetilde{A}_s := A_s - a_1c_1$ and $G_2 = \widetilde{A}_t$ for some $s,t \geq 3$ such that $n = s+ t -1$.
    It then follows from \eqref{eq:addingtheodds} that $A_{n} \cong G_1 \cup G_2$.
    By the induction hypothesis, we have $c_2(G_1 + a_1c_1)$ and $c_2(G_2 + a_1c_1)$ are odd.
    Furthermore, $c_2(G_1 + b_1c_1)$ and $c_2(G_2 + b_1c_1)$ are even by \eqref{eq:Rmult}. 
    Hence we can apply \Cref{thm:oddglue} to obtain that either
    \[
    H = G_1 \cup G_2 \cong A_{n} \; \text{ or } \; \widetilde{H} = G_1 \cup \widetilde{G}_2
    \]
    has odd 2-realisation number.
    However, \Cref{thm:oddglue} gives us no indication which of the two graphs has odd 2-realisation number, and we cannot deduce that it is $A_n$ without the additional analysis of \Cref{sec:infiniteoddgraphs}.
\end{remark}

\subsection*{Acknowledgements}

We thank Matthias Adrian-Himmelmann for assistance with computations.
S.\,D.\ was partially supported by the Heilbronn Institute for Mathematical Research,
and supported by the KU Leuven grant iBOF/23/064 and the FWO grants G0F5921N (Odysseus) and G023721N.
A.\,N.\ and B.\,S.\ were partially supported by EPSRC grant EP/X036723/1. A.\,N.\ was partially supported by UK Research and Innovation (grant number UKRI1112), under the EPSRC Mathematical Sciences Small Grant scheme.

\bibliographystyle{plainurl}
{\small
\bibliography{ref}
}

\end{document}